\documentclass[10pt,USletter,twocolumn]{article}
\usepackage[utf8]{inputenc}

\usepackage{orcidlink}
\usepackage{setspace}

\usepackage[electronic]{ifsym}
\usepackage{indentfirst,latexsym}
\usepackage{bm}
\usepackage{mathrsfs}
\usepackage{amsmath,amsfonts,amsthm,amscd,amssymb}
\usepackage{url}

\usepackage{graphicx}
\usepackage{subfigure}

\usepackage[all,pdf]{xy}

\usepackage{booktabs}
\usepackage{multirow}
\usepackage{multicol}
\usepackage{stfloats}  
\usepackage{color}

\usepackage{algorithm}         
\usepackage{algorithmicx}
\usepackage{algpseudocode}
\usepackage{float}
\usepackage{bm}
\usepackage{bbm}

\newcommand{\Theo}{\textbf{Theorem}~}
\newcommand{\Algr}{\textbf{Algorithm}~}
\newcommand{\Fig}{\textbf{Figure}~}
\newcommand{\Tab}{\textbf{Table}~}
\newcommand{\cpvar}[1]{\texttt{#1}}

\newcommand{\ProcName}[1]{\textsc{#1}}
\newcommand{\Alg}{\textsc{Alg}}

\algnewcommand\algorithmicswitch{\textbf{switch}}
\algnewcommand\algorithmiccase{\textbf{case}}
\algnewcommand\algorithmicassert{\texttt{assert}}
\algnewcommand\Assert[1]{\State \algorithmicassert(#1)}%
\algdef{SE}[SWITCH]{Switch}{EndSwitch}[1]{\algorithmicswitch\ #1\ \algorithmicdo}{\algorithmicend\ \algorithmicswitch}%
\algdef{SE}[CASE]{Case}{EndCase}[1]{\algorithmiccase\ #1}{\algorithmicend\ \algorithmiccase}%
\algtext*{EndCase}%

\makeatletter
\renewcommand{\ALG@name}{Algorithm}
\newenvironment{breakablealgorithm}
{
	\begin{center}
		\refstepcounter{algorithm}
		\setlength{\baselineskip}{15pt} 
		\renewcommand{\caption}[2][\relax]{
			\hrule height.9pt depth0pt \kern3pt
			{\raggedright\textbf{\ALG@name~\thealgorithm} ##2\par}%
			\ifx\relax##1\relax 
			\addcontentsline{loa}{algorithm}{
				\protect\numberline{\thealgorithm}##2}%
			\else 
			\addcontentsline{loa}{algorithm}{
				\protect\numberline{\thealgorithm}##1}%
			\fi
			\kern2pt\hrule\kern2pt
		}
	}{
		\kern3pt\hrule\relax
	\end{center}
}

\DeclareMathOperator{\dif}{d}  
\DeclareMathOperator{\rank}{rank}  
\DeclareMathOperator{\Ker}{Ker}

\newcommand{\scrd}[2]{{#1}_{\mathrm{#2}}}

\renewcommand{\vec}[1]{\bm{#1}}

\newcommand{\mpair}[2]{ \left\langle {#1}, {#2} \right\rangle}

\newcommand{\mat}[1]{\bm{#1}}

\newcommand{\set}[1]{\left\{ #1 \right\}}
\newcommand{\seq}[1]{\langle #1 \rangle}
\newcommand{\abs}[1]{\left| #1 \right|}

\newcommand{\norm}[1]{\left\lVert #1 \right\rVert}
\newcommand{\normp}[2]{{\left\lVert #1 \right\rVert}_{#2}}
\newcommand{\trsp}[1]{#1^\textsf{T}}

\newcommand{\inv}[1]{#1^{-1}}
\newcommand{\ginv}[1]{#1^\dagger}    
\newcommand{\tinv}[1]{#1^{-\textsf{T}}}

\newcommand{\st}{\mathrm{s.t.} }
\newcommand{\sign}{\textrm{sign} }

\newcommand{\GL}[2]{\mathrm{GL}(#1,\mathbb{#2})}               
\newcommand{\ES}[3]{\mathbbm{#1}^{{#2}\times {#3}}}     

\newcommand{\mybmatrix}[1]{\left[#1\right]}

\newcommand{\DRL}{\mathrm{DRL}}
\newcommand{\Prox}{\mathrm{Prox}}
\newcommand{\indcc}[1]{\mathcal{I}_{#1}}

\newcommand{\blue}[1]{{#1}}

\newtheorem{thm}{Theorem}
\renewcommand{\Theo}{\textbf{Theorem} }

\usepackage[left=1.5cm,right=1.5cm,top=2cm,bottom=2cm]{geometry}

\title{A Unified Theoretic and Algorithmic Framework  for Solving  Multivariate Linear Model with $\ell^1$-norm Approximation\thanks{Corresponding author: Hong-Yan Zhang,\\ e-mail: hongyan@hainnu.edu.cn}}

\author{Zhi-Qiang Feng$^a$\orcidlink{0000-0002-7384-295X}, Hong-Yan Zhang$^{*a}$\orcidlink{0000-0002-4400-9133}, Ji Ma$^{b}$\orcidlink{0009-0009-3345-7430}, \\
Daniel Delahaye$^c$\orcidlink{0000-0002-4965-6815}, Ruo-Shi Yang$^d$\orcidlink{0009-0003-0210-9309} and Man Liang$^e$\orcidlink{0000-0003-0577-0832}\\
\begin{tabular}{l}
$^a$\small{\textit{School of Information Science and Technology, Hainan Normal University, Haikou 571158, China}}\\
$^b$\small{\textit{Sino-european Institute of Aviation Engineering (SIAE), Civil Aviation University of China, Tianjin 300300, China}}\\
$^c$\small{\textit{Lab ENAC, École Nationale de l’Aviation Civile (ENAC), Toulouse 31400, France}}\\
$^d$\small{\textit{College of Civil Aviation, Nanjing University of Aeronautics and Astronautics,  Nanjing 210016, China}}\\
$^e$\small{\textit{Department of Aerospace Engineering and Aviation, RMIT University, Melbourne, VIC 3000, Australia}}\\
\end{tabular}
}

\date{May 20, 2025}

\begin{document}
\maketitle 

\begin{abstract}
It is a challenging problem that solving the \textit{multivariate linear model} (MLM) $\mat{A}\vec{x}=\vec{b}$ with the $\ell_1 $-norm approximation method such that $\normp{\mat{A}\vec{x}-\vec{b}}{1}$, the $ \ell_1 $-norm of the \textit{residual error vector} (REV), is minimized. In this work, our contributions lie in two aspects: firstly, the equivalence theorem for the structure of the $ \ell_1 $-norm optimal solution to the MLM is proposed and proved; secondly, a unified algorithmic framework for solving the MLM with $\ell_1$-norm optimization is proposed and six novel algorithms (L1-GPRS, L1-TNIPM, L1-HP, L1-IST, L1-ADM, L1-POB) are designed. There are three significant characteristics in the algorithms discussed: they are implemented with simple matrix operations which do not depend on specific  optimization solvers; 
they are described with algorithmic pseudo-codes and implemented with Python and Octave/MATLAB which means easy usage; and the high accuracy and efficiency of our six new algorithms can be achieved successfully in the scenarios with different levels of data redundancy. We hope that the unified theoretic and algorithmic framework with source code released on GitHub could motivate the applications of the $\ell_1$-norm optimization for parameter estimation of MLM arising in science, technology, engineering, mathematics, economics, and so on. 
\end{abstract}

\tableofcontents

\section{Introduction}
\label{sec:introduction}
It is a key issue that how to solve the \textit{multivariate linear model} (MLM)
accurately, efficiently and robustly in science, technology, economics, management, and so on \cite{McCullagh1989-2ed,Dobson2018-4ed}. The MLM is also named by overdetermined system of linear algebraic equations  and  \textit{generalized linear model} (GLM). Formally, the MLM can be expressed by
\begin{equation}\label{eq-MLM}
	\mat{A}\vec{x}=\vec{b}
\end{equation}
where $\vec{x} = (x_i)_{n\times 1}\in \ES{R}{n}{1}$ is the unknown $n$-dim parameter vector, $\mat{A} = (a_{ij})_{m\times n}\in \ES{R}{m}{n}$ is the coefficient matrix such that $m \ge n$, and $\vec{b} = (b_i)_{m\times 1}\in \ES{R}{m}{1}$ is the $m$-dim observation vector.

Although the form of MLM \eqref{eq-MLM} arising in various applications is simple, it is not well-defined \cite{CADZOW2002524} due to the noise existing in the matrix $\mat{A}$ and/or in the vector $\vec{b}$. 
\footnote{Conceptually, the MLM \eqref{eq-MLM} is called well-defined or consistent if there exists at least one solution $ \vec{x} $ such that the equation \eqref{eq-MLM} holds on strictly in the sense of algebra. However, it is not our topic in this study.}
Usually, the vector
\begin{equation}\label{intro-3}
	 \vec{r(x)}=\mat{A}\vec{x} -\vec{b}\in \ES{R}{m}{1}
\end{equation}
is called the \textit{residual error vector} (REV). The non-negative function
\begin{equation}
\mathcal{C}_p(\vec{x}) = \normp{\vec{r}(\vec{x})}{p}^p, \quad \vec{x}\in \ES{R}{n}{1}, p\ge 1
\end{equation} 
is called the cost function for the MLM, in which $\normp{\cdot}{p}$ is the 
$\ell^p$-norm defined by 
\begin{equation}\label{intro-2}
	\normp{\vec{w}}{p}=\sqrt[p]{\sum_{i = 1}^{m}\abs{w_i}^p}, \quad \forall\  \vec{w}\in\ES{R}{m}{1}.
\end{equation}
In order to reduce the impact of the noise arising in $\mat{A}$ and/or $\vec{b}$, it is necessary to solve the optimization problem
  \begin{equation} \label{eq-mlm-opt}
  \scrd{\vec{x}}{opt} = \arg \min_{\vec{x}\in \ES{R}{m}{1}} \mathcal{C}_p(\vec{x}) 
  \end{equation}
with the given integer $p\ge 1$. The most popular choice  for the cost function is 
\begin{equation} \label{eq-CostFun-l2}
\mathcal{C}_2(\vec{x}) = \normp{\mat{A}\vec{x}-\vec{b}}{2}^2 = \trsp{(\mat{A}\vec{x}-\vec{b})}(\mat{A}\vec{x}-\vec{b}),
\end{equation}
which is a smooth function.  This choice leads to the family of least square methods, which includes the classic \textit{least squares} (LS) \cite{GolubLoan2013,ZhangXD2017Matrix} and its variants such as \textit{data least squares} (DLS), \textit{total least squares} and \textit{scaled total least squares} (STLS) \cite{Paige02a:STLS-Fundmentals,Paige02b:STLS-Bounds,Paige02c:STLS-unifying-LS-TLS-DLS,Zhang2009novel,Zhang2022multi}. 
The advantages of the family of least square methods include the following aspects: 
\begin{itemize}
\item firstly the differentiability of $\mathcal{C}_2(\vec{x})$ in \eqref{eq-CostFun-l2} leads to a closed form of optimal solution to \eqref{eq-mlm-opt} derived by the orthogonality principle in Hilbert space; 
\item secondly, the  geometric interpretations for the TLS  is intuitive, and 
\item finally the statistical and topological interpretation for the STLS are clear.
\end{itemize}
However, there are still some disadvantages for the family of least square methods: 
\begin{itemize}
\item[i)] when the distribution of the noise in $\mat{A}$ and/or  $\vec{b}$ is not the normal distribution, the performance of the estimation can not be guaranteed;
\item[ii)] if there are some outliers, the parameter estimation obtained by the family of least square methods will be useless.  
\end{itemize}

In order to avoid the impact of outliers, the RANSAC method was proposed by Fischler and Bolles in 1981 \cite{Fischler1981RANSAC}, which is widely used in computer vision. Generally, the time computational complexity of RANSAC method and its variations is high due to the process of random sampling if the ratio of outliers is high.
For the purpose of seeking precise and robust solution  to  \eqref{eq-mlm-opt} effectively by dealing with the dark sides i) and ii) simultaneously, it is a good choice to take the $\ell^1$-norm  instead of the $\ell^2$-norm. In this case,  the cost function 
\begin{equation}
\mathcal{C}_1(\vec{x}) = \normp{\mat{A}\vec{x}-\vec{b}}{1} = \sum_{i=1}^{m}\abs{\sum_{j=1}^{n}a_{ij}x_j-b_i}
\end{equation}
is not differentiable, which leads to the following challenging optimization issue
\begin{equation}\label{P_1}\tag{$ \rm{P}_1 $}
	\scrd{\vec{x}}{opt}=\min_{\vec{x} \in \ES{R}{n}{1}} \mathcal{C}_1(\vec{x}).
\end{equation}
The solution to the problem \eqref{P_1} is known as solving the MLM  in the sense of $\ell^1$-norm minimization. The solution is called the \textit{minimum $\ell^1$-norm approximate solution} to the MLM. Simultaneously, the problem \eqref{P_1} is called the \textit{minimum $\ell^1$-norm approximation problem} \cite{CADZOW2002524}.

Wagner \cite{wagner1959linear} in 1959 highlighted the correlation between linear programming and $ \ell^1 $-norm approximation problem. In 1966, Barrodale et al. \cite{barrodale1966algorithms} introduced a primal algorithm that capitalizes on the unique structure of the linear programming formulation for \eqref{P_1}, and later presented an enhanced version in \cite{Barrodale1973}. Converting the $ \ell^1 $-norm approximation problem into a linear programming problem offers a refined mathematical representation, albeit at the expense of increasing the complexity problem-solving (see section \ref{sec:L1andLP}).

In 1967, Usow \cite{usow1967l_1} introduced a descent method for calculating the optimal $ \ell^1 $-norm approximation. This method involves locating the lowest vertex of a convex polytope that represents the set of all possible $ \ell^1 $-norm approximations. In 2002, Cadzow et al. \cite{CADZOW2002524} proposed an improved $ \ell^1 $-norm perturbation algorithm based on brute-force methods for solving for the minimum $ \ell^1 $-norm approximate solution (see section \ref{sec:L1PTB}). Although the structures for their iterative algorithms are similar, these algorithms are limited due to the high computational complexity when the dimension $n$ is large.

Bartels et al. in 1978 presented a projection descent method by minimizing piecewise differentiable cost functions. This method achieved the minimization of $ \mathcal{C}_1(\vec{x}) $ in a finite number of steps. In 1991, Yang and Xue \cite{yzh1991-en} proposed an active set algorithm  by converting the \eqref{P_1} into a piecewise linear programming problem via establishing an active equation index set and solving it  iteratively. In 1993, Wang and Shen \cite{wjs1993-en} proposed an interval algorithm to solve the minimum $ \ell^1 $-norm approximation problem. However, the structure of the interval algorithm is complex and the existence of the optimal solution interval should be verified repeated. Unfortunately, the interval algorithm is not attractive until now.  

Yao \cite{yjk2007-en} proposed an effective numerical method  for problem \eqref{P_1} based on the minimum $\ell^1$-norm residual vector in 2007. By converting the problem \eqref{P_1} to a constrained $ \ell^1 $-norm optimization problem, Yao obtained the minimum $ \ell^1 $-norm approximate solution through a compatible system of linear equations with the assumption $\rank(\mat{A})=n$. However, Yao's method does not hold when $\rank(\mat{A}) <n$. 

\begin{figure*}[htbp]
    \centering
    \includegraphics[width=0.95\textwidth]{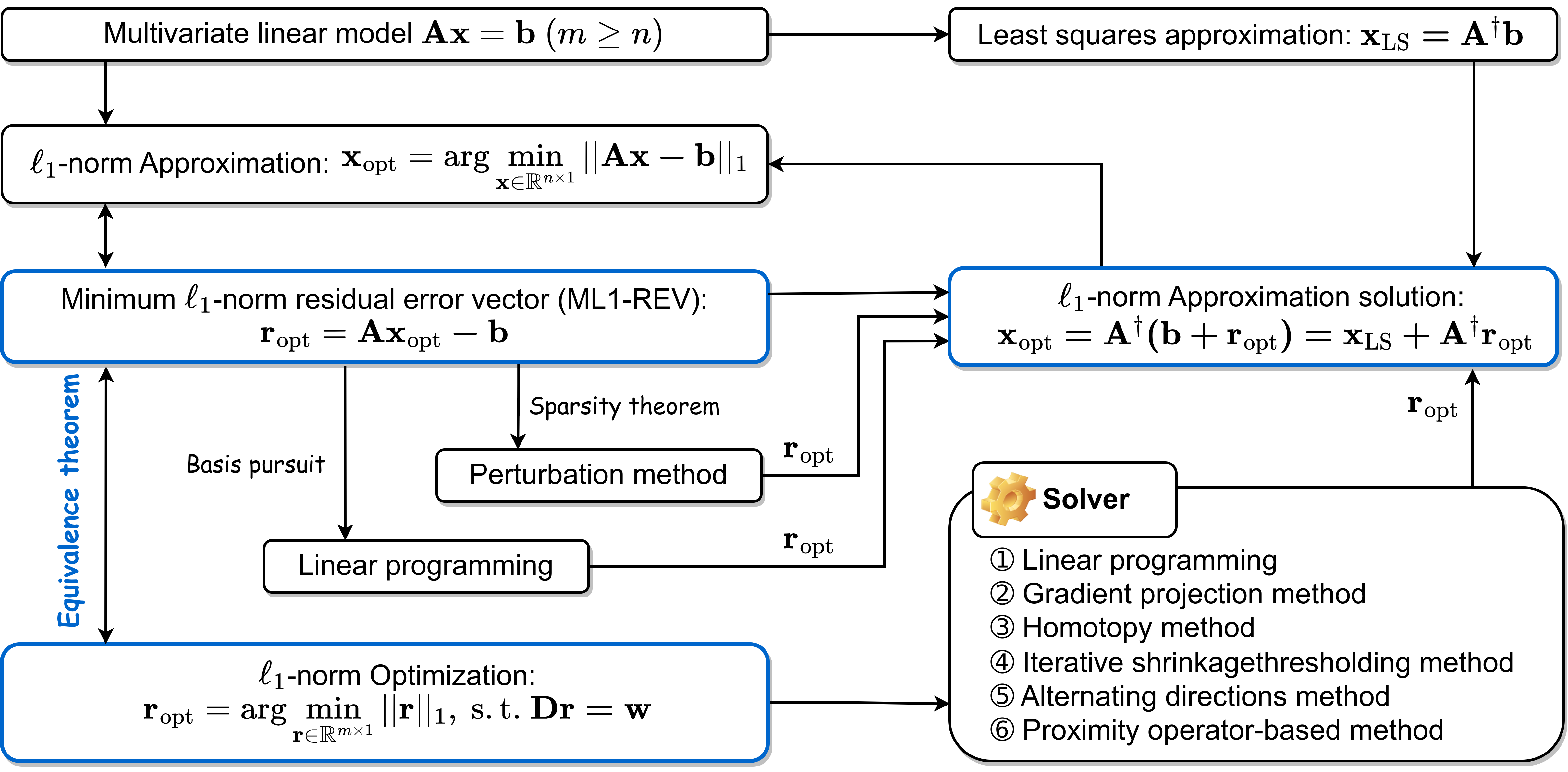}
    \caption{Unified theoretic and algorithmic framework for solving the MLM with $\ell_1$-norm optimization} 
    \label{fig-flowchart}
\end{figure*}

In summary, there is a lack of good method and ready-to-use algorithms for solving the $\ell^1$-norm optimization problem \eqref{P_1} with the general condition $\rank(\mat{A}) \le n$, in which the "good" means high precision, high robustness,  and low or acceptable computational complexity. 
In this study, our contributions for solving \eqref{P_1} include the following aspects: 
\begin{itemize}
\item[i)] a unified theoretic framework is given by establishing and proving the equivalence theorem, which states the new sufficient and necessary condition for solving the minimum $ \ell^1 $-norm approximate solution with the basis pursuit method and the  Moore-Penrose inverse; 
\item[ii)] a unified algorithmic framework is proposed via REV, which covers all of the available algorithms for \eqref{P_1}; 
\item[iii)] six new algorithms are designed, which just relies on fundamental matrix operations and do not depend on specific optimization solvers; 
\item[iv)]  the performance of different algorithms are evaluated with numerical simulations; and
\item[v)] all of the algorithms designed are described with pseudo-code in the sense of computer science and the Python/Octave implementations are released on the GitHub, which not only reduces the difficulty of understanding the mathematical principles but also increases the potential usability and popularity in the sense of engineering applications. 
\end{itemize}
The global view of this work is illustrated in \Fig\ref{fig-flowchart}, which includes the principle and algorithms for solving the MLM.

The subsequent contents of this study are organized as follows: Section \ref{sec-preliminaries} introduces the preliminaries for this paper; Section \ref{sec-prop-ell1opt} copes with the properties of $\ell^1$-norm approximate solution to the MLM; Section \ref{sec-framework} deals with the unified algorithmic framework for solving the \eqref{eq-mlm-opt} via REV; Section \ref{sec-alg-performance} concerns the performance evaluation for the algorithms proposed; and finally Section \ref{sec-conclusion} gives the conclusions. 

For the convenience of reading, the notations and corresponding interpretations are summarized in \Tab \ref{table-1}.
\begin{table*}[htb]
\centering
\tabcolsep=3pt
\renewcommand{\baselinestretch}{1.2}
\caption{Mathematical notations}\label{table-1}
\resizebox{\textwidth}{!}{
\begin{tabular}{cl}
\toprule
\textbf{Symbol}  & \textbf{Interpretation} \\ \midrule
$u^+$  & positive part of $u$: $\forall u\in \mathbb{R}, u^+=\max(u,0)$ \\
$u^-$  & negative part of $u$: $\forall u\in \mathbb{R}, u^-=\max(-u,0)$ \\
$u^+ - u^-$ & decomposition of $u\in \mathbb{R}$ with positive and negative parts such that $u = u^+ -u^-$\\
$\abs{u}$ & absolution of $u$: $\forall u\in \mathbb{R}, \abs{u} = u^+ + u^-$ \\
$\vec{v}^+=\max(\vec{v},0)$ & positive part of $\vec{v} = (v_i)_{m\times 1}\in \ES{R}{m}{1}$ such that $v^+_i = \max(v_i,0)$ for $1\le i\le m$\\
$\vec{v}^-=\max(-\vec{v},0)$ & negative part of $\vec{v} = (v_i)_{m\times 1}\in \ES{R}{m}{1}$ such that $v^-_i = \max(-v_i,0)$ for $1\le i\le m$\\
$\vec{v}^+ - \vec{v}^-$ & decomposition of $\vec{v}\in \ES{R}{m}{1}$ with positive and negative parts such that $\vec{v} = \vec{v}^+ -\vec{v}^-$\\
$\vec{1}_m$ & $m$-dim column vector with unit components: $\vec{1}_m = \trsp{[1, 1, \cdots, 1]}\in \ES{R}{m}{1}$ \\ 
$ \ginv{\mat{A}} $ & Moore-Penrose inverse of matrix $ \mat{A} $\\
$ \mat{I}_n $ & Identity matrix of order $n $\\
$ \vec{e}_k $ & Standard basis vector whose components are all 0 except for its $ k $-th component which is 1\\
$\mathscr{Z} $ & index set $\mathscr{Z}=\set{i_1, \cdots, i_k, \cdots, i_{m_0}}$\\
$\mathscr{Z}^c $ & complementary set of $\mathscr{Z}$, i.e., $\mathscr{Z}^c = \set{1,2,\cdots, m}-\mathscr{Z}  = \set{j_1, \cdots, j_t, \cdots, j_{m-m_0}}$ such that $ m_0\le m$\\
$ \vec{a}(i\!:\!j) $ & sub-vector $\vec{a}(i\!:\!j) = \trsp{[a_i, a_{i+1}, \cdots, a_j]}$ constructed from the components of vector $\vec{a}\in \ES{R}{n}{1}$ such that $1\le i < j \le n$ \\
$ \vec{a}\succcurlyeq 0 $ & non-negative vector $\vec{a}=\trsp{[a_1, a_2, \dots, a_n]}$ such that $a_i\ge 0$ for $1\le i\le n$\\
$ \mat{A}(i_1\!\!:\!i_2,j_1\!\!:\!j_2) $ & sub-block of $\mat{A}\in \ES{R}{m}{n}$ specified by the row indices $i_1\sim i_2$ and column indices $j_1\sim j_2$ such that $1\le i_1 < i_2 \le m, 1\le j_1 < j_2 \le n$\\
$ \mat{A}(i_1\!\!:\!i_2,:) $ & sub-block of $\mat{A}\in \ES{R}{m}{n}$ specified by all rows whose indices range from $i_1\sim i_2$, where $ 1\le i_1< i_2\le m $\\
$ \mat{A}(:,j_1\!\!:\!j_2) $ & sub-block of $\mat{A}\in \ES{R}{m}{n}$ specified by all columns whose indices range from $j_1\sim j_2$, where $1\le j_1 < j_2 \le n$\\
$\mat{A}_z=\mat{A}[\mathscr{Z}]$ & sub-block of $\mat{A}\in \ES{R}{m}{n}$ specified by the row index $\mathscr{Z}$, i.e., $ \mat{A}[\mathscr{Z}] = \mat{A}(\mathscr{Z}, :) = \mat{A}(\mathscr{Z}, 1:n)$ \\
$\mat{A}_*=\mat{A}[\mathscr{Z}^c]$ & sub-block of $\mat{A}\in \ES{R}{m}{n}$ specified by the row index $\mathscr{Z}^c$, i.e., $ \mat{A}[\mathscr{Z}^c] = \mat{A}(\mathscr{Z}^c, :) = \mat{A}(\mathscr{Z}^c, 1:n)$ \\
$\mat{A} =\mpair{\mat{A}_z}{\mat{A}_*}$ & decomposition of $\mat{A}\in \ES{R}{m}{n}$ such that $\mat{A}_*=\mat{A}[\mathscr{Z}]$ and $\mat{A}_z=\mat{A}[\mathscr{Z}^c]$ \\
$\vec{b} =\mpair{\vec{b}_z}{\vec{b}_*}$ & decomposition of $\vec{b}\in \ES{R}{m}{1}$ such that $\vec{b}_z = \vec{b}[\mathscr{Z}]$ and $\vec{b}_z=\vec{b}[\mathscr{Z}^c]$  \\
$\vec{r} =\mpair{\vec{r}_z}{\vec{r}_*}$ & decomposition of REV $\vec{r}\in \ES{R}{m}{1}$ such that $\vec{r}_z = \vec{r}[\mathscr{Z}]$ and $\vec{r}_z=\vec{r}[\mathscr{Z}^c]$ \\
$ \vec{a}\odot\vec{b} $&the element-wise product for two vectors $ \vec{a} $ and $ \vec{b} $ of the same dimension $ n\times1 $, with elements given by $ (\vec{a}\odot\vec{b})_i=a_i b_i,1\le i\le n $\\
$ \vec{a}\oslash\vec{b} $&the element-wise division for two vectors $ \vec{a} $ and $ \vec{b} $ of the same dimension $ n\times1 $, with elements given by $ (\vec{a}\oslash\vec{b})_i=\dfrac{a_i}{b_i}(b_i\neq0),1\le i\le n $\\
\bottomrule
\end{tabular}
}
\end{table*}

\section{Preliminaries} \label{sec-preliminaries}

\subsection{Notations}

For any real number $u\in \mathbb{R}$, it can be decomposed into its positive part $u^+$ and negative part $u^-$. Formally, we have
\begin{equation}\label{notation-1}
u = u^+ - u^-, \quad \abs{u} = u^+ + u^-
\end{equation}
where
\begin{equation*}\label{notation-2}
u^+ = \max(u, 0), \quad u^- = \max(-u, 0). 
\end{equation*}
 As an extension of \eqref{notation-1}, for any positive integer $m\in \mathbb{N}$ and any real vector $\vec{v}=(v_i)_{m\times 1}\in \ES{R}{m}{1}$, we have
\begin{equation}\label{notation-3}
\vec{v} = \vec{v}^+ - \vec{v}^-  = (v^+_i)_{m\times 1} - (v^-_i)_{m\times 1}
\end{equation}
where
\begin{equation}\label{notation-4}
\vec{v}^+ = \begin{bmatrix}
v^+_1 \\
v^+_2 \\
\vdots\\
v^+_m
\end{bmatrix} \succcurlyeq \vec{0}, \quad
\vec{v}^-= \begin{bmatrix}
v^-_1 \\
v^-_2 \\
\vdots \\
v^-_m
\end{bmatrix} \succcurlyeq \vec{0}
\end{equation}
are the positive and negative parts of $\vec{v}$. Let
\begin{equation}\label{notation-5}
\vec{1}_m=\trsp{[1,1,\cdots,1]}\in\ES{R}{m}{1}
\end{equation}
be the vector with $m$ components of $1$, then  \eqref{notation-3} and \eqref{notation-4} imply that
\begin{equation}\label{notation-6}
\normp{\vec{v}}{1} = \sum^m_{i=1}\abs{v_i} = \sum^m_{i=1}(v^+_i + v^-_i) = \trsp{\vec{1}}_m(\vec{v}^+ + \vec{v}^-).
\end{equation}
For the parameter $\vec{x}$ and the observation vector $\vec{b}$ in \eqref{eq-MLM}, equation \eqref{notation-3} shows that
\begin{equation}\label{notation-7}
		\vec{x} = \vec{x}^+ - \vec{x}^-, \quad 
		\vec{r} = \vec{r}^+ - \vec{r}^-.
\end{equation}
With the help of \eqref{intro-3} we can obtain
\begin{equation}\label{notation-8}
	\mat{A}(\vec{x}^+ - \vec{x}^-)  - (\vec{r}^+ - \vec{r}^-)= \vec{b}.
\end{equation}

%

\subsection{Matrix Decomposition Induced by REV}

\begin{figure*}[htb]
	\centering
	\includegraphics[width=.8\linewidth]{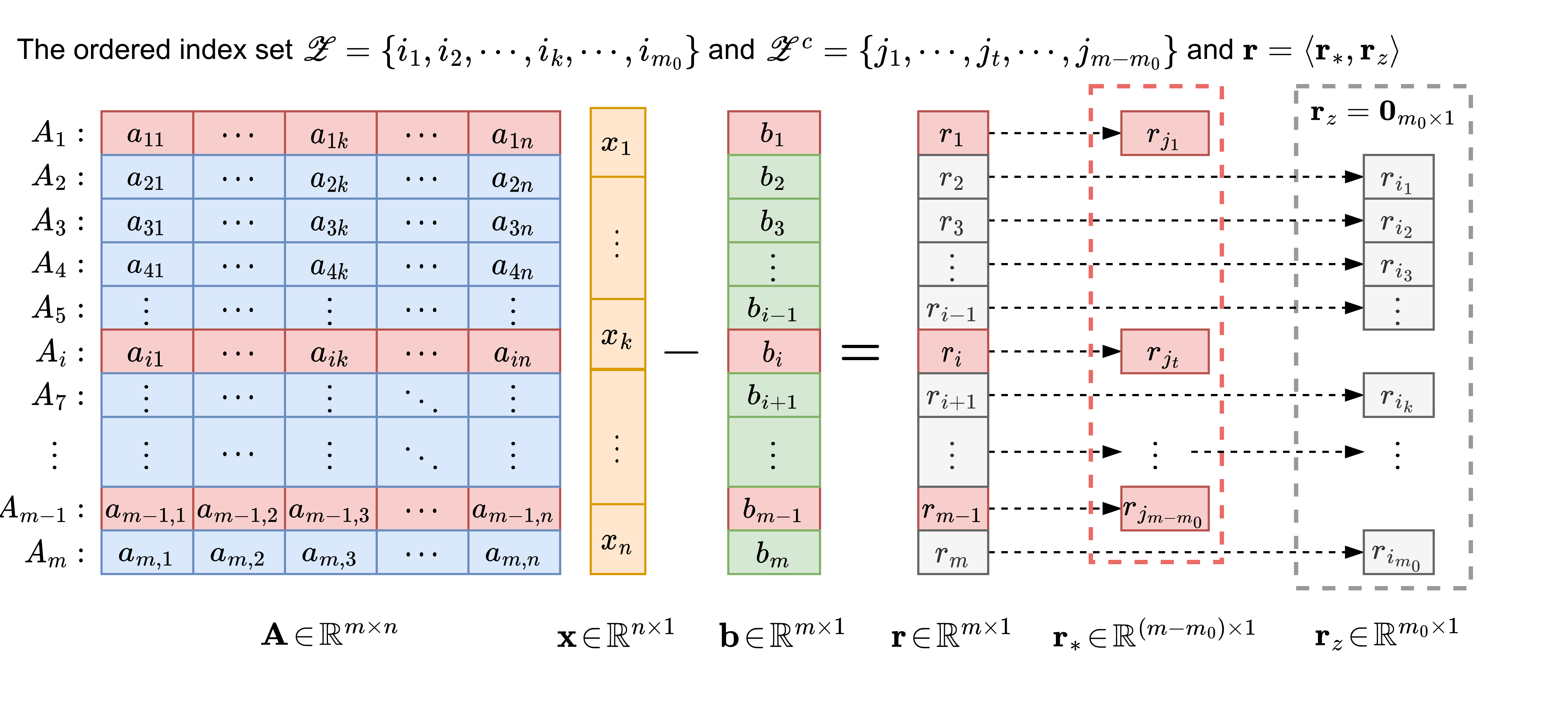}	
	\caption{Decompositon of the REV $ \vec{r} = \mat{A}\vec{x} -\vec{b} = \mpair{\vec{r}_*}{\vec{r}_z} = \mpair{\vec{r}[\mathscr{Z}]}{\vec{r}[\mathscr{Z}^c]}$} \label{sprase_decop}
\end{figure*}

Suppose that there are $ m_0\in\set{0,1,\cdots,m} $ zeros in the components of REV $\vec{r} = (r_i)_{m\times 1} = \mat{A}\vec{x} -\vec{b}\in \ES{R}{m}{1}$, i.e.
\begin{equation}
	r_{i_k} = 0, \quad 1\le i_1<i_2<\cdots < i_{m_0} \le m.
\end{equation}
If $m_0 = m$ then $\vec{r}=\vec{0}$ and $\mat{A}\vec{x}=\vec{b}$ must be compatible. 
Let 
\begin{equation}
\mathscr{Z} = \set{t: r_t = 0, 1\le t\le m} =\set{i_1, \cdots, i_k, \cdots, i_{m_0}}
\end{equation}
be an ordered index set and 
\begin{equation}
\begin{aligned}
\mathscr{Z}^c  &= \set{t: r_t \neq 0, 1\le t\le m}\\
&=\set{j_1,\cdots,j_t, \cdots,j_{m-m_0}}
\end{aligned}
\end{equation}
be the complementary set of $\mathscr{Z}$. The matrix $\mat{A}$ can be decomposed by
\begin{equation}
\mat{A} = \mpair{\mat{A}_z}{\mat{A}_*} 
\end{equation}
such that
\begin{equation}\label{decomp:A_z}
\begin{cases}
\mat{A}_z = \mat{A}[\mathscr{Z}]=\mat{A}(\mathscr{Z},1:n)\in\ES{R}{m_0}{n}\\
\mat{A}_* = \mat{A}[\mathscr{Z}^c]=\mat{A}(\mathscr{Z}^c,1:n)\in\ES{R}{(m-m_0)}{n}.
\end{cases}
\end{equation}
In other words, $\mat{A}_z$ consists of the $i_1$-th, $i_2$-th, $\cdots$, $i_{m_0}$-th columns of $\mat{A}$. Similarly, the vectors $\vec{b}$ and $\vec{r}$ can be decompose by
\begin{equation}\label{decomp-b}
\left\{
\begin{split}
\vec{b} &= \mpair{\vec{b}_z}{\vec{b}_*} = \mpair{\vec{b}[\mathscr{Z}]}{\vec{b}[\mathscr{Z}^c]}\\
\vec{r} &= \mpair{\vec{r}_z}{\vec{r}_*} = \mpair{\vec{r}[\mathscr{Z}]}{\vec{r}[\mathscr{Z}^c]}
\end{split}
\right.
\end{equation}
such that
\begin{equation} \label{decomp-z}
\begin{cases}
\vec{r}_z(\vec{x}) = \mat{A}_z\vec{x} - \vec{b}_z = \vec{0}\in \ES{R}{m_0}{1}\\
\vec{r}_*(\vec{x}) = \mat{A}_* \vec{x} -\vec{b}_* \neq \vec{0}\in \ES{R}{(m-m_0)}{1}
\end{cases}
\end{equation}
where $ \vec{b}_{z} $ and $ \mat{A}_{z} $ correspond to the subset of $ m_0 $ equations which have zero residuals ,  $ \vec{b}_* $ and $ \mat{A}_* $ correspond to the entities of the remaining subset of $ m-m_0 $ equations which have non-vanishing residuals. \Fig \ref{sprase_decop} shows the decomposition of the REV intuitively by separating the vanishing part $\vec{r}_z$ in which each component is zero and the non-vanishing part $\vec{r}_*$ in which each component is not zero. When the value $(m- m_0)/m$ is small, the REV is called \textit{sparse REV}.

\subsection{Minimum $ \ell^1 $-norm REV}

\subsubsection{Concept} 

Suppose that $ \scrd{\vec{x}}{opt} $ is the solution to \eqref{P_1}, the vector
\begin{equation}
\scrd{\vec{r}}{opt} = \mat{A}\scrd{\vec{x}}{opt} - \vec{b} 
\end{equation}
is called the \textit{minimum $ \ell^1 $-norm residual error vector} (ML1-REV) of \eqref{P_1}. This definition was introduced by Cui \& Quan \cite{cmg1996-en} in 1996. Obviously, if the ML1-REV $ \scrd{\vec{r}}{opt}$ is known, then the solution to the consistent system of linear equations $ \mat{A}\vec{x} = \vec{b} + \scrd{\vec{r}}{opt} $ will be the solution to \eqref{P_1}.

\subsubsection{Equivalence and Sparsity}

The properties of equivalence and sparsity for the MLM by Feng \& Zhang  \cite{Fzq2022} and Cadzow \cite{CADZOW2002524} are proved respectively, which are stated by \Theo \ref{fzq-theorem-2} and \Theo \ref{cadzow-theorem-1} respectively. 

\begin{thm}[Equivalence]\label{fzq-theorem-2}
	Let $ \mat{Ax = b} $ and $ \mat{Bx = b} $ be two different MLM. If $ \exists \mat{P}\in \GL{n}{R}=\set{\mat{Q}\in \ES{R}{n}{n}, \det(\mat{Q})\neq 0}$ such that $ \mat{B = AP}$, then the ML1-REV of the two MLMs are the same.
\end{thm}

\begin{thm}[Sparsity]\label{cadzow-theorem-1}
	For any $ \vec{b}\in\ES{R}{m}{1} $ and any $ \mat{A}\in\ES{R}{m}{n} $ such that $m\ge n$, there exists $ \vec{x}_0 \in \ES{R}{n}{1}$ which minimizes the cost function $ \mathcal{C}_1(\vec{x}) $ such that the REV
	\begin{equation*}\label{perturbation-1}
		\vec{r}(\vec{x}_0)=\mat{A}\vec{x}_0-\vec{b}
	\end{equation*}
has at least $ n $ zero components. Furthermore, if the row vectors of the augmented matrix $ \mybmatrix{\mat{A},\vec{b}}\in\ES{R}{m}{(n+1)} $ satisfy the Haar condition then there exists $ \vec{x}_0\in \ES{R}{n}{1} $ which minimizes $ \mathcal{C}_1(\vec{x}) $ and the $\vec{r}(\vec{x}_0)$ has exactly $n$ zero components.
\end{thm}

\subsubsection{Perturbation Strategy}\label{sec:L1PTB}

Cadzow \cite{CADZOW2002524} proposed the perturbation strategy by constructing an iterative method for solving the \eqref{P_1}.

\begin{thm}[Perturbation Strategy]\label{cadzow-theorem-2}
For the MLM $\mat{A}\vec{x}=\vec{b}$, suppose that
\begin{itemize}
\item[i)] there are $ m_0 $ zeroes indicated by $\mathscr{Z}$ in the REV $ \vec{r}(\vec{x}) = \mat{A}\vec{x}-\vec{b}$ where $ 0 \le m_0 < n $;
\item[ii)] the decompositions induced by $\mathscr{Z}$ are $\mat{A}=\mpair{\mat{A}_z}{\mat{A}_*}$, $\vec{b}=\mpair{\vec{b}_z}{\vec{b}_*}$ and $\vec{r} = \mpair{\vec{r}_z}{\vec{r}_*}$  respectively.
\end{itemize}
For the direction vector  
	\begin{equation}\label{perturbation-3}
	\vec{d}\in \Ker(\mat{A}_z) =\set{\vec{p}\in \ES{R}{n}{1}:	\mat{A}_z\vec{p}=\vec{0}}
	\end{equation}
and the optimal step size
	\begin{equation}\label{perturbation-4}
		\alpha = \arg\min_{\gamma_k \in W}\normp{\vec{r}_*(\vec{x})+\gamma_k\mat{A}_*\vec{d}}{1}
	\end{equation}
where
	\begin{equation} \label{perturbation-alpha-k}
	W=\set{\gamma_k=-\dfrac{\trsp{\vec{e}_k}\vec{r}_*(\vec{x})}{\trsp{\vec{e}}_k\mat{A}_*\vec{d}}: 1\le k\le m, \trsp{\vec{e}_k}\mat{A}_*\vec{d}\neq 0}
	\end{equation}
and $\vec{e}_k$ is the $k$-th column of the $m\times m$ identity matrix, the direction $\alpha\vec{d}$ reduces the cost function $\mathcal{C}_1(\vec{x})$  by  
\begin{equation}
\mathcal{C}_1(\vec{x} + \alpha \vec{d}) \le \mathcal{C}_1(\vec{x})
\end{equation}
and the REV $ \vec{r}(\vec{x}+\alpha\vec{d}) $ has at least $ m_0 + 1 $ zero components.
\end{thm}

\Theo \ref{cadzow-theorem-1} indicates that the optimal $\vec{x}$ which minimizing $\mathcal{C}_1(\vec{x})$ ensures that $\abs{\mathscr{Z}}\ge n $, viz., the number of the zeros in REV is at least $n$ for the MLM restricted by $m\ge n$. 

\Theo \ref{cadzow-theorem-2} presents an iterative strategy with finite times of perturbation via the feasible decreasing direction $\alpha \vec{d}$ in each iteration.  The iteration must stop if $ m_0=n $, thus it is necessary to add a step in order to determine whether the current choice of $ \vec{x} $ is the minimum $ \ell^1 $-norm solution of interest by constructing an alternative feasible direction for the iteration.  Bloomfield \& Steiger \cite{bloomfield1980least} proposed the following three steps strategy for this purpose:
\begin{itemize}
\item firstly, constructing the auxiliary decision vector
\begin{equation}\label{perturbation-5}
	\vec{s}=\tinv{\mat{A}_z}\trsp{\mat{A}_*}\cdot\sign\set{\mat{A}_*\vec{x}-\vec{b}_{*}}\in\ES{R}{m_0}{1}
\end{equation}
which results the uniqueness of the optimal solution if $ \normp{\vec{s}}{\infty}=\max\limits_{i}\abs{s_i}\le 1 $;
\item secondly, computing the new feasible direction by
\begin{equation}\label{perturbation-6}
	\alpha\vec{d}=c\inv{\mat{A}_z}\vec{u}(\vec{s}),\quad c\in\mathbbm{R}
\end{equation}
where the vector $\vec{u}$ depends on the vector $\vec{s}$ such that its components can be specified by
\begin{equation}\label{perturbation-7}
	u_i(\vec{s})=\left\{\begin{aligned}
		&1,&\abs{s_i}>1,\\
		&0,&\abs{s_i}\le1.
	\end{aligned}\right.,\quad 1 \le i \le m_0
\end{equation}
\item thirdly, updating the $\vec{x}$ with $\vec{x} + \alpha \vec{d}$ according to  \eqref{perturbation-6} after 
$\vec{s}$ and $\vec{d}$ being computed from \eqref{perturbation-5} and \eqref{perturbation-7}. 
\end{itemize}


The perturbation strategy based on Cadzow's \Theo \ref{cadzow-theorem-2} for solving the $\ell^1$-norm approximation solution to the MLM is presented in \Algr \ref{alg-Perturbation}. The \verb|CBS| in the procedure name \ProcName{L1ApproxPertCBS} comes from the names of Cadow, Bloomfield and Steiger.
\begin{breakablealgorithm}
	\caption{Solving  the $ \ell^1 $-approximation via perturbation.}\label{alg-Perturbation}
	\begin{algorithmic}[1]
		\Require $ \mat{A}\in\ES{R}{m}{n},\vec{b}\in\ES{R}{m}{1},c>0$ and $m>n\ge 2 $, $ \cpvar{maxiter}\in\mathbb{N}^{+} $.
		\Ensure $ \vec{x}\in\ES{R}{n}{1} $.
		\Function{L1ApproxPertCBS}{$ \mat{A}, \vec{b},c $,\cpvar{maxiter}}
		\State Initialize $ m $ as the number of rows of matrix $ \mat{A} $, $ n $ as the number of columns;
		\State $ \vec{x}\gets\vec{0}_n $;\Comment{initialize}
        \State $ \cpvar{iter}\gets0 $;
		\While{$ \cpvar{iter} < \cpvar{maxiter} $}
        \State $ \cpvar{iter}\gets\cpvar{iter}+1 $;
		\State $\vec{r}\gets\mat{A}\vec{x}-\vec{b}$;
		\State $ \mathscr{Z}\gets\set{i: 1\le i\le m, r_i=0} $;
		\While{$ \abs{\mathscr{Z}}<n $}
		\State $\mat{A}_z \gets \mat{A}[\mathscr{Z}]$;
		\State $\mat{A}_* \gets \mat{A}[\mathscr{Z}^c]$;
		\State $\vec{r}_* \gets \vec{r}[\mathscr{Z}^c]$;
		\State Compute the kernel $\Ker(\mat{A}_z)$;
		\State Choose a vector $ \vec{d}\in\Ker(\mat{A}_z) $;\Comment{See \eqref{perturbation-3}}
		\State $ \vec{v},\vec{f}\gets\vec{0}_{\abs{\mathscr{Z}^c}} $;\Comment{initialize with zeros}
		\For{$ i\in\set{1,2,\cdots,\abs{\mathscr{Z}^c}} $}
		\State $ \vec{e}\gets\vec{0}_{\abs{\mathscr{Z}^c}} $;
		\State $ e_{i}\gets1 $;\Comment{$ \forall\ \trsp{\vec{e}}_k\mat{A}_*\vec{d}\neq0 $}
		\If{$ \trsp{\vec{e}}\mat{A}_*\vec{d}==0 $}
		\State $ f_i\gets+\infty $;
		\Else
		\State $ v_i\gets-\trsp{\vec{e}}\vec{r}_*/\trsp{\vec{e}}\mat{A}_*\vec{d} $;
		\State $ f_i\gets\normp{\vec{r}_*+v_i\mat{A}_*\vec{d}}{1} $;
		\EndIf
		\EndFor
		\State $ \set{f_{\min},i_{\min}}\gets\ProcName{SearchMin}(\vec{f}) $;\Comment{See \eqref{perturbation-4}}
		\State $ \vec{x}\gets\vec{x}+v_{i_{\min}}\cdot\vec{d} $;
		\State $\vec{r}\gets\mat{A}\vec{x}-\vec{b}$;
		\State $ \mathscr{Z}\gets\set{i: 1\le i\le m, r_i=0}$;
		\EndWhile
		\State $\mat{A}_z \gets \mat{A}[\mathscr{Z}]$;
		\State $\mat{A}_* \gets \mat{A}[\mathscr{Z}^c]$;
		\State $\vec{r}_* \gets \vec{r}[\mathscr{Z}^c]$;
		\State $\vec{s}\gets\tinv{\mat{A}_z}\trsp{\mat{A}_*}\cdot\sign(\vec{r}_*) $;\Comment{Get $ \vec{s} $ via \eqref{perturbation-5}}
		\If{$ \normp{\vec{s}}{\infty}\le1 $}
		\State \textbf{break};
		\Else
		\State $ \vec{e}^{(s)}\gets\vec{0}_{\abs{\mathscr{Z}}} $;\Comment{Get $ \vec{e}^{(s)} $ from \eqref{perturbation-7}}
		\For{$ i\in\set{1,2,\cdots,\abs{\mathscr{Z}}} $}
		\If{$ \abs{s_i}>1 $}
		\State $ e^{(s)}_{i}\gets1 $;
		\EndIf
		\EndFor
		\State $ \vec{x}\gets\vec{x}+c\inv{\mat{A}_z}\vec{e}^{(s)} $;\Comment{Update $\vec{x}$ by \eqref{perturbation-6}}
		\EndIf
		\EndWhile
		\State\Return $\vec{x} $;
		\EndFunction
	\end{algorithmic}
\end{breakablealgorithm}

\subsection{Linear Programming and MLM }\label{sec:L1andLP}

Barrodale and Roberts in \cite{Barrodale1973,Barrodale1974}  reformulated the problem \eqref{P_1} as the following standard linear programming problem
\begin{equation}\label{LP_1}\tag{$ \rm{LP}_{1} $}
	\begin{aligned}
		&\min\limits_{\vec{r}^{+},\vec{r}^{-}\in\ES{R}{m}{1}}\trsp{\vec{1}}_m(\vec{r}^{+}+\vec{r}^{-})\\
		&\st \left\{\begin{aligned}
			&\mat{A}(\vec{x}^{+}-\vec{x}^{-})-(\vec{r}^{+}-\vec{r}^{-})=\vec{b}\\
		&\vec{x}^{+}\succcurlyeq 0\\
		&\vec{x}^{-}\succcurlyeq 0\\
		&\vec{r}^{+}\succcurlyeq 0\\
		&\vec{r}^{-}\succcurlyeq 0
		\end{aligned}\right.
	\end{aligned}
\end{equation}
by splitting the vectors $\vec{r}, \vec{x}, \vec{b}$ with the help of \eqref{notation-6}, \eqref{notation-7} and \eqref{notation-8}. The alternative form of \eqref{LP_1} can be expressed by
\begin{equation} \label{linprog-2}
\min_{\vec{y}\in \ES{R}{(2m+2n)}{1}} \trsp{\vec{c}}\vec{y}, \quad \st \quad \scrd{\mat{A}}{eq}\vec{y}=\vec{b}, ~\vec{y}\succcurlyeq \vec{0}
\end{equation}
where
\begin{equation}
\scrd{\mat{A}}{eq} = [-\mat{I}_m, \mat{I}_m, \mat{A}, -\mat{A}]\in \ES{R}{m}{(2m+2n)}
\end{equation}
and
\begin{equation}
\vec{c} = \begin{bmatrix}
\vec{1}_m \\ 
\vec{1}_m\\
\vec{0}_n\\
\vec{0}_n
\end{bmatrix}, \quad 
\vec{y} = \begin{bmatrix}
\vec{r}^{+}\\ \vec{r}^{-}\\ \vec{x}^{+}\\ \vec{x}^{-}
\end{bmatrix}\in \ES{R}{(2m+2n)}{1}.
\end{equation}
by reformulating the matrices and vectors involved. 
Thus the available optimization algorithms such as the interior-point method or simplex method \cite{Barrodale1973} can be utilized to solve \eqref{LP_1}.

The procedure \ProcName{L1ApproxLinProg} listed in \Algr \ref{alg-L1linprog} is specifically designed to solve the $ \ell^1 $-approximation solution using a standard linear programming solver.
\begin{breakablealgorithm}
	\caption{Solving the $\mat{A}\vec{x}=\vec{b}$ with $ \ell^1 $-norm approximation via linear programming.}\label{alg-L1linprog}
	\begin{algorithmic}[1]
		\Require $ \mat{A}\in\ES{R}{m}{n},\vec{b}\in\ES{R}{m}{1}$ and $m>n\ge 2 $. 
		\Ensure $ \scrd{\vec{x}}{opt}\in\ES{R}{n}{1} $ in the sense of $\ell^1$-norm optimization.
		\Function{L1ApproxLinProg}{$ \mat{A}, \vec{b} $}
		\State $[m,n]\gets \cpvar{size}(\mat{A})$;\Comment{the size of the matrix $\mat{A}$}
		\State $d\gets 2m+2n$;\Comment{the dimension of $\vec{y}$}
		\State $ \trsp{\vec{c}} \gets\mybmatrix{\trsp{\vec{1}}_m,\trsp{\vec{1}}_m,\trsp{\vec{0}}_n,\trsp{\vec{0}}_n}\in \ES{R}{1}{d}$;
		\State $ \scrd{\mat{A}}{eq}\gets [-\mat{I}_m, \mat{I}_m, \mat{A}, -\mat{A}]\in \ES{R}{m}{d} $;
		\State $ \scrd{\vec{y}}{opt}\!\gets\! \ProcName{LinProgSolver}(\trsp{\vec{c}}, \scrd{\mat{A}}{eq}, \vec{b}, \vec{0}_{d\times 1}) $;
		\State $ \scrd{\vec{x}}{opt} \gets \scrd{\vec{y}}{opt}(2m+1:2m+n)-\scrd{\vec{y}}{opt}(2m+n+1:d) $;
		\State\Return $ \scrd{\vec{x}}{opt} $;
		\EndFunction
	\end{algorithmic}
\end{breakablealgorithm}

The \ProcName{LinProgSolver} in \Algr \ref{alg-L1linprog} is used to solve problem \eqref{linprog-2}, which can be obtained from lots of available standard tools and packages\footnote{The open-source version in Python is available at \url{https://docs.scipy.org/doc/scipy/reference/generated/scipy.optimize.linprog.html}.}.

\section{Theoretic Framework of $\ell^1$-norm Approximation Solution to MLM}  \label{sec-prop-ell1opt}
 
\subsection{Structure of $\ell^1$-norm Approximate Solution}

Given the optimal REV $\vec{r} = \vec{r}(\vec{x})$ for \eqref{P_1} in the sense of $\ell^1$-norm optimization, it is necessary to find the source $\vec{x}$ from the image $\vec{r}$. Our exploration shows that the optimal solution to \eqref{P_1} can be obtained through computing the Moore-Penrose inverse of the matrix $\mat{A}$. Actually, we have the following theorem for the structure of the $\ell^1$-norm approximate solution.
\begin{thm}  \label{fzq-theorem-1}
	Suppose $m, n \in \mathbb{N}$ and $m \ge n\ge 2$, for $\vec{x} \in \mathbb{R}^{n \times 1}$, $\vec{b} \in \mathbb{R}^{m \times 1}$ and  $\mat{A} = \begin{bmatrix}		\mat{A}_1\\ \mat{A}_2 \end{bmatrix} \in \mathbb{R}^{m \times n}$ such that $\mat{A}_1\in \ES{R}{n}{n}$ and $\mat{A}_2\in \ES{R}{(m-n)}{n}$.  Let
	\begin{equation}\label{residual-1}
		\mat{D} = \begin{bmatrix} 
		-\mat{A}_2\ginv{\mat{A}_1} & \mat{I}_{m-n}
		\end{bmatrix} \in \ES{R}{(m-n)}{n}
	\end{equation}
and
	\begin{equation}\label{residual-2}
		\vec{w} = \mat{A}_2\mat{A}_1^\dagger\vec{b}(1:n) - \vec{b}(n+1:m) \in \mathbb{R}^{(m-n) \times 1}
	\end{equation}
where $\ginv{(\cdot)}$ is the Moore-Penrose inverse and $\mat{I}_{m-n}$ is the identity matrix of order $m-n$. 
	Solving the minimum $\ell^1$-approximation solution to the MLM $\mat{A}\vec{x}=\vec{b}$ is equivalent to solving the ML1-REV
	\begin{equation}\label{BP}\tag{BP}
		\scrd{\vec{r}}{opt} = \arg \min_{\vec{r}\in \ES{R}{m}{1}}\norm{\vec{r}}_1,\ \st \ \mat{D}\vec{r}=\vec{w}
	\end{equation}
	and we have
	\begin{equation}\label{residual-4}
		\scrd{\vec{x}}{opt} =\mat{{A}}^\dagger(\vec{b}+\scrd{\vec{r}}{opt}). 
	\end{equation}
\end{thm}


\proof 

The definition of Moore-Penrose inverse implies that 
\begin{equation} \label{eq-A2-A1-ginv}
\mat{A}_2=\mat{A}_2\mat{I}_{n} = \mat{A}_2\ginv{\mat{A}_1}\mat{A}_1.
\end{equation} 
Substituting \eqref{eq-A2-A1-ginv} into the block form of $\mat{A}$, we immediately have 
\begin{equation} \label{eq-A-new-form}
\mat{A} = \begin{bmatrix}
\mat{A}_1 \\ \mat{A}_2
\end{bmatrix}
 =\begin{bmatrix}
\mat{A}_1 \\ \mat{A}_2\ginv{\mat{A}_1}\mat{A}_1
\end{bmatrix} = \begin{bmatrix}
\mat{I}_n \\ \mat{A}_2\ginv{\mat{A}_1} 
\end{bmatrix} \mat{A}_1
\end{equation}
Hence
\begin{equation}\label{proof-1}
	{\begin{bmatrix}
			\mat{I}_n\\\mat{A}_2\ginv{\mat{A}_1}
	\end{bmatrix}}\mat{A}_1\vec{x}=\vec{b}.
\end{equation}
by \eqref{eq-A-new-form} and \eqref{eq-MLM}.
For the ML1-REV 
\begin{equation*}\label{proof-2}
	\vec{r}=\trsp{\mybmatrix{r_1, r_2, \cdots, r_m}} 
\end{equation*}
of the noisy MLM, \Theo \ref{fzq-theorem-2} shows that \eqref{proof-1} shares the same ML1-REV with the MLM \eqref{eq-MLM}. In other words, $\vec{r}$ is also the ML1-REV of \eqref{proof-1}. Put 
\begin{equation}
\mat{A}_1\vec{x}=\vec{z},
\end{equation} 
then the compatible linear system of equations for the MLM \eqref{eq-MLM} can be written by
\begin{equation}\label{proof-3}
	\begin{bmatrix}
		\mat{I}_n\\\mat{A}_2\ginv{\mat{A}_1}
	\end{bmatrix}\vec{z}=\vec{b} + \vec{r}.
\end{equation}
Splitting the vectors $\vec{b}$ and $\vec{r}$ of length $m$ into two sub-vectors of length $n$ and $m-n$ respectively,  \eqref{proof-3} can be written by
\begin{equation}\label{proof-4}
	\begin{bmatrix}
		\vec{z} \\ \mat{A}_2\ginv{\mat{A}_1}\vec{z}
	\end{bmatrix}
	=\begin{bmatrix}
		\vec{b}(1:n)+\vec{r}(1:n)\\
		\vec{b}(n+1:m) + \vec{r}(n+1:m)
	\end{bmatrix}.
\end{equation}
The substitution of $\vec{z}$ with \eqref{proof-4} yields
\begin{equation}\label{proof-6}
	\begin{aligned}
		&\mat{A}_2\ginv{\mat{A}_1}(\vec{b}(1:n)+\vec{r}(1:n))
		=\vec{b}(n+1:m)+\vec{r}(n+1:m).
	\end{aligned}
\end{equation}
Let
\begin{equation}\label{proof-7}
	\mat{C}=\mat{A}_2\mat{A}_1^\dagger\in\ES{R}{(m-n)}{n}
\end{equation}
and substitute \eqref{residual-2} into \eqref{proof-6}, we can obtain 
\begin{equation}\label{proof-9}
	\begin{aligned}
		\vec{w} = \vec{r}(n+1:m)-\mat{C}\vec{r}(1:n)
	\end{aligned}
\end{equation}
or equivalently 
\begin{equation}\label{proof-10}
	r_{n+i}-\sum_{j=1}^{n}C_{ij}r_j=w_i,\quad i=1,2,...,m-n.
\end{equation}
Consequently, \eqref{proof-10} and \eqref{proof-9} are equivalent to the linear constraint
\begin{equation}
	\mat{D}\vec{r} = \vec{w}
\end{equation}
according to \eqref{residual-1}. Q.E.D.

The equation \eqref{residual-4} in \Theo \ref{fzq-theorem-1} implies that the $\ell^1$-norm approximation solution $\scrd{\vec{x}}{opt}$ to \eqref{eq-MLM} consists of the traditional LS solution $\ginv{\mat{A}}\vec{b}$ and the correction term $\ginv{\mat{A}}\scrd{\vec{r}}{opt}$ specified by the ML1-REV.

\subsection{MLM and Penalized Least Squares Problem}

It is evident that the $ \ell^1 $-approximation problem \eqref{P_1} has been transformed into the standard \textit{basis pursuit} (BP) \cite{chen2001atomic} problem \eqref{BP}. For the problem \eqref{BP}, we can relax the equality constraint as follows
\begin{equation}\label{BPDN}\tag{$ \rm{BP}_\epsilon $}
	\scrd{\vec{r}}{opt} = \arg \min_{\vec{r}\in \ES{R}{m}{1}}\norm{\vec{r}}_1,\ \st \ \normp{\mat{D}\vec{r}-\vec{w}}{2}\le\epsilon
\end{equation}
where $ \epsilon\ge0 $. Using a Lagrangian formulation, the problem \eqref{BPDN} can be reformulated into a penalized least squares problem, viz.
\begin{equation}\label{QP_lambda}\tag{$ \rm{QP}_\lambda $}
	\scrd{\vec{r}}{opt} = \arg \min_{\vec{r}\in \ES{R}{m}{1}} \frac{1}{2}\normp{\mat{D}\vec{r}-\vec{w}}{2}^2+\lambda\norm{\vec{r}}_1
\end{equation}
where $ \lambda\ge0 $ is the Lagrangian multiplier. By setting $ \epsilon\rightarrow0 $ and $ \lambda\rightarrow0 $, it becomes clear that the solutions of problem \eqref{BPDN} and \eqref{QP_lambda} coincide with those of the problem \eqref{BP} \cite{yang2010review}. The extensive and well-established research about the \eqref{BPDN} and \eqref{QP_lambda} can be used to deal with the problem \eqref{BP} as described in \cite{chen2015fast}.

\section{Algorithmic Framework of $ \ell^1 $-norm Approximation Solution to MLM}
\label{sec-framework}

\subsection{Engineering of Available $\ell^1 $-norm Optimization Methods for Solving the REV}

Currently, there are at least six methods for solving the REV in the sense of $\ell^1 $-norm optimization. Usually, these methods are presented in a mathematical style instead of engineering style since the algorithmic pseudo-codes are missing. With the purpose of reducing the difficulty of applying the $\ell^1 $-norm optimization in complex engineering problems, we give brief description about the mathematical principles and present clear algorithmic pseudo-codes for the methods in various literature. 

We remark that the objective of engineering education is to train the students' ability of solving complex engineering and technical problems with the conceive-design-implement-operate (CDIO) approach \cite{CdioSyllabus2022V3,HYZhang2024CAE}, in which "design" involves algorithm design via algorithmic pseudo-codes.

\subsubsection{Linear Programming Method}

Feng et al. transformed the problem \eqref{BP} into the standard linear programming problem in order to solve the minimum $ \ell^1 $-norm REV \cite{Fzq2022}. The procedure \ProcName{L1OptLinProg} listed in \Algr \ref{alg-BPlinprog} is designed to solve the minimum $ \ell^1 $-norm REV by calling the standard linear programming solver with the interface \ProcName{LinProgSolver}($\vec{t}$, \textrm{Aeq}, \textrm{beq}, \textrm{lb}). The mathematical principles are presented in Appendix \ref{sec:LinprogL1Opt} with details.
\begin{breakablealgorithm}
	\caption{Solving the $\ell^1 $-norm REV via linear programming.}\label{alg-BPlinprog}
	\begin{algorithmic}[1]
		\Require $ \mat{D}\in\ES{R}{(m-n)}{m},\vec{w}\in\ES{R}{(m-n)}{1} $.
		\Ensure $ \vec{r}\in\ES{R}{m}{1} $.
		\Function{L1OptLinProg}{$ \mat{D}, \vec{w} $}
		\State $ \vec{t} \gets\trsp{\vec{1}}_{2m}$;
		\State $ \mat{\Phi}\gets [\mat{D}, -\mat{D}]\in \ES{R}{(m-n)}{2m} $;
		\State $ \scrd{\vec{\beta}}{opt}\gets$\ProcName{LinProgSolver}($\vec{t}$, $\mat{\Phi}$, $\vec{w}$, $\vec{0}_{2m\times 1}$); 
		\State $ \scrd{\vec{r}}{opt} \gets \scrd{\vec{\beta}}{opt}(1:m)-\scrd{\vec{\beta}}{opt}(m+1:2m) $;\Comment{by \eqref{BP-linprog-4}}
		\State\Return $ \scrd{\vec{r}}{opt} $;
		\EndFunction
	\end{algorithmic}
\end{breakablealgorithm}

\subsubsection{Gradient Projection Method}

Two gradient descent methods can be used to solve problem \eqref{QP_lambda} and obtain the REV, namely the \textit{Gradient Projection Sparse Representation} (GPSR) method \cite{figueiredo2007gradient} and the \textit{Truncated Newton Interior-point Method} (TNIPM) \cite{kim2007interior}. The implementations of the GPSR method and TNIPM have been released in a open way by the authors\footnote{A MATLAB implementation of GPSR is available at \url{http://www.lx.it.pt/~mtf/GPSR/}. A MATLAB Toolbox for TNIPM, called L1LS, is available at \url{https://web.stanford.edu/~boyd/l1_ls/}.}. We provide two standard function interfaces for the GPSR method and TNIPM for solving  \eqref{QP_lambda}:
\begin{itemize}
	\item $ \ProcName{L1OptGPSR}(\mat{D},\vec{w},\lambda) $, which solves the problem \eqref{QP_lambda} with the GPSR method described in \cite{figueiredo2007gradient};
	\item $ \ProcName{L1OptTNIPM}(\mat{D},\vec{w},\lambda) $, which solves the problem \eqref{QP_lambda} with the Preconditioned Conjugate Gradients(PCG) accelerated version of the TNIPM\cite{kim2007interior}.
\end{itemize}
The corresponding pseudo-code for the GPSR method is provided in \Algr \ref{alg-GPSRBB}. More details about the mathematical principles for both two methods are described in Appendix \ref{sec:GPL1Opt}.

\begin{breakablealgorithm}
	\caption{Solving the $\ell^1$-norm REV via the GPSR method with Barzilai-Borwein gradient projection.}\label{alg-GPSRBB}
	\begin{algorithmic}[1]
		\Require $ \mat{D}\in\ES{R}{(m-n)}{m},\vec{w}\in\ES{R}{(m-n)}{1},\lambda>0,\varepsilon>0 $, $ \cpvar{maxiter}\in\mathbb{N}^{+} $.
		\Ensure $ \vec{r}\in\ES{R}{m}{1} $.
		\Function{L1OptGPSR}{$ \mat{D}, \vec{w}, \lambda, \varepsilon $, \cpvar{maxiter}}
        \State $\alpha \gets 1.0$;\Comment{set parameters}
        \State $\alpha_{\min} \gets 10^{-30}$;
        \State $\alpha_{\max} \gets 10^{30}$;
        \State $\vec{r} \gets \vec{0}_{m\times 1}$;\Comment{initialization}
        \State $\vec{u} \gets  \vec{0}_{m\times 1}$; \Comment{$\vec{u}$ is for $\vec{r}^+$}
        \State $\vec{v} \gets  \vec{0}_{m\times 1}$; \Comment{$\vec{v}$ is for $\vec{r}^-$}
        \State $\vec{\mu} \gets \mat{D} \vec{r}$;
        \State $ \cpvar{iter}\gets0 $;
        \While{$\cpvar{iter}<\cpvar{maxiter}$}
        \State $ \cpvar{iter}\gets\cpvar{iter}+1 $;
        \State $\nabla_{\vec{u}}Q \gets \trsp{\vec{D}}(\mat{D} \vec{r} - \vec{w}) + \lambda$;\Comment{compute gradients}
        \State $\nabla_{\vec{v}}Q \gets -\nabla_{\vec{u}}Q + 2\lambda$; \Comment{$-\trsp{\vec{D}}(\mat{D} \vec{r} - \vec{w}) + \lambda$}        
        \State $\dif\vec{u} \gets (\vec{u} - \alpha \cdot \nabla_{\vec{u}}Q)^+ - \vec{u}$;\Comment{search direction}
        \State $\dif\vec{v} \gets (\vec{v} - \alpha \cdot \nabla_{\vec{v}}Q)^+ - \vec{v}$;\Comment{search direction}
        \State $\dif\vec{r} \gets \dif\vec{u} - \dif\vec{v}$;
        
        \State $\gamma \gets \normp{\mat{D} \cdot \dif\vec{r}}{2}^2$;
        \State $ \beta_0\gets-\dfrac{1}{\gamma}{\left[\trsp{(\nabla_{\vec{u}}Q)} \cdot \dif\vec{u} + \trsp{(\nabla_{\vec{v}}Q)} \cdot \dif\vec{v}\right]} $;\Comment{$ \gamma\neq0 $}
        \State $\beta \gets \min\left(\beta_0, 1\right)$;
        
        \State $\vec{u}_{\text{improve}} \gets \vec{u} + \beta \cdot \dif\vec{u}$;\Comment{update parameters}
        \State $\vec{v}_{\text{improve}} \gets \vec{v} + \beta \cdot \dif\vec{v}$;
        \State $\vec{y} \gets \min\left(\vec{u}_{\text{improve}}, \vec{v}_{\text{improve}}\right)$;
        \State $\vec{u} \gets \vec{u}_{\text{improve}} - \vec{y}$;
        \State $\vec{v} \gets \vec{v}_{\text{improve}} - \vec{y}$;
        \State $\vec{r} \gets \vec{u} - \vec{v}$;
        
        \State $\delta \gets (\dif\vec{u})^T \cdot \dif\vec{u} + (\dif\vec{v})^T \cdot \dif\vec{v}$;
        \If{$\gamma \leq 0$}\Comment{compute new alpha}
        \State $\alpha \gets \alpha_{\text{max}}$;
        \Else
        \State $\alpha \gets \min\left(\alpha_{\text{max}}, \max\left(\alpha_{\min}, \dfrac{\delta}{\gamma}\right)\right)$;
        \EndIf
        \State $\vec{\mu} \gets \vec{\mu} + \beta \cdot \mat{D} \cdot \dif\vec{r}$;
        \If{$\dfrac{\normp{\dif\vec{r}}{2}}{\normp{\vec{r}}{2}} \leq \varepsilon$}
        \State \textbf{break}
        \EndIf
        \EndWhile

		\State\Return $ \vec{r} $;
		\EndFunction
	\end{algorithmic}
\end{breakablealgorithm}

\begin{breakablealgorithm}
	\caption{Solving the $\ell^1$-norm REV via the TNIPM with Preconditioned Conjugate Gradients algorithm.}\label{alg-TNIPM}
	\begin{algorithmic}[1]
		\Require $ \mat{D} \in \mathbb{R}^{(m-n) \times m}, \vec{w} \in \mathbb{R}^{(m-n)\times 1}, \lambda>0, \varepsilon>0 $, $ \cpvar{maxiter}\in\mathbb{N}^{+} $.
		\Ensure $ \vec{r} \in \mathbb{R}^{m\times 1} $.
		\Function{L1OptTNIPM}{$ \mat{D}, \vec{w}, \lambda, \varepsilon $, $ \cpvar{maxiter} $}
		\State $ \mu \gets 2 $;\Comment{Initialize parameters}
		\State $ \alpha \gets 0.01 $;
		\State $ \beta \gets 0.5 $;
		\State $ \vec{r} \gets \vec{0}_{m\times1} $;
		\State $ \vec{u} \gets \vec{1}_{m\times1} $;
		\State $ \vec{f}\gets \begin{bmatrix}
		\vec{u}-\vec{r}\\\vec{u}+\vec{r}
		\end{bmatrix} $;
		\State $ \dif\vec{f} \gets \vec{0}_{2m\times1} $;
		
		\State $ s \gets +\infty $;
		\State $ d_{\cpvar{obj}} \gets -\infty $;
		\State $ t \gets \min\left(\max\left(1, \dfrac{1}{\lambda}\right), \dfrac{2m}{\varepsilon}\right) $;\Comment{parameter for the primal interior-point method}
		
		\State $ \cpvar{iter}\gets0 $;
        \While{$\cpvar{iter}<\cpvar{maxiter}$}
        \State $ \cpvar{iter}\gets\cpvar{iter}+1 $;
		\State $ \vec{z} \gets \mat{D} \vec{r} - \vec{w} $;
		\State $ \vec{\nu} \gets \vec{z} $;\Comment{Construct a dual feasible point}
		\If{$ \normp{\trsp{\mat{D}} \vec{\nu}}{\infty} > \lambda $}
		\State $ \vec{\nu} \gets \vec{\nu} \cdot \frac{\lambda}{\normp{\trsp{\mat{D}} \vec{\nu}}{\infty}} $;
		\EndIf
		\State $ p_{\cpvar{obj}} \gets 0.5 \trsp{\vec{z}} \vec{z} + \lambda \norm{\vec{r}}_1 $;\Comment{primary objective}
		\State $ d_{\cpvar{obj}} \gets \max(-0.5\trsp{\vec{\nu}} \vec{\nu} - \trsp{\vec{\nu}} \vec{w}, d_{\cpvar{obj}}) $;\Comment{dual objective}
		\State $ \eta \gets p_{\cpvar{obj}} - d_{\cpvar{obj}} $;\Comment{duality gap}
		
		\If{$ \dfrac{\eta}{d_{\cpvar{obj}}} < \varepsilon $}\Comment{Check stopping criterion}
		\State \Return $ \vec{r} $;
		\EndIf
		
		\If{$ s \geq 0.5 $}\Comment{Update t}
		\State $ t \gets \max \left( \mu \cdot \min \left( \dfrac{2m}{\eta}, t \right), t \right) $;
		\EndIf
		
		\State $ \vec{q}_1 \gets {\vec{1}}\oslash{(\vec{u} + \vec{r})} $;
		\State $ \vec{q}_2 \gets {\vec{1}}\oslash{(\vec{u} - \vec{r})} $;
		\State $ \nabla F \gets \begin{bmatrix}
		\trsp{\mat{D}} \vec{z} - \dfrac{\vec{q}_1 - \vec{q}_2}{t}\\ \lambda \cdot \vec{1} - \dfrac{\vec{q}_1 + \vec{q}_2}{t}
		\end{bmatrix} $;\Comment{gradient}
		\State $ \mat{B}_1 \gets \dfrac{1}{t}\cdot\operatorname{diag}(\vec{q}_1\odot\vec{q}_1 + \vec{q}_2\odot\vec{q}_2) $; 
		\State $ \mat{B}_2 \gets \dfrac{1}{t}\cdot\operatorname{diag}(\vec{q}_1\odot\vec{q}_1 - \vec{q}_2\odot\vec{q}_2) $;
		\State $ \mat{H}\gets\begin{bmatrix}
		\trsp{\mat{D}}\mat{D}+\mat{B}_1&\mat{B}_2\\
		\mat{B}_2&\mat{B}_1
		\end{bmatrix} $;\Comment{hessian matrix}
		\State $ \mat{P}\gets\begin{bmatrix}\mat{I}_{m}+\mat{B}_1&\mat{B}_2\\ \mat{B}_2&\mat{B}_1\end{bmatrix} $;\Comment{Compute preconditioner}
		\State $ \cpvar{tol} \gets \min\left(0.1, \dfrac{\varepsilon \cdot \eta}{\min\left(1, \normp{\nabla F}{2}\right)}\right) $;\Comment{Set preconditioned conjugate gradient tolerance}
		
		\State $ \dif\vec{f} \gets \ProcName{PCG}(\mat{H},-\nabla F,\mat{P},\dif\vec{f},\cpvar{tol}) $;\Comment{Solve the system of linear equations $ \mat{P}^{-1}\mat{H}\dif\vec{f}=-\mat{P}^{-1}\nabla F $ for the optimal search direction $ \dif\vec{f} $ using PCG algorithm with preconditioner $ \mat{P} $ and initial guess $ \dif\vec{f} $;}
		
		\State $ \dif\vec{r} \gets \dif\vec{f}(1:m)$;\Comment{Split results}
		\State $ \dif\vec{u} \gets \dif\vec{f}(m+1:2m) $;
		
		\State $  F \gets 0.5 \vec{z}^T \vec{z} + \lambda \sum\limits_{i=1}^{m}u_i - \sum\limits_{i=1}^{2m} \dfrac{\log(f_i)}{t} $;
		\State $ s \gets 1.0 $;\Comment{the step size}
		\While{\cpvar{TRUE}}\Comment{Backtracking line search}
		\State $ \vec{r}' \gets \vec{r} + s \cdot \dif\vec{r} $;
		\State $ \vec{u}' \gets \vec{u} + s \cdot \dif\vec{u} $;
		\State $ \vec{f}' \gets \begin{bmatrix}
		\vec{u}' - \vec{r}'\\ \vec{u}' + \vec{r}'
		\end{bmatrix} $;
		\If{$ \max(\vec{f}') > 0 $}
		\State $ \vec{z}' \gets \mat{D} \vec{r}' - \vec{w} $;
		\State $  F' \gets 0.5 (\vec{z}')^T \vec{z}' + \lambda \sum\limits_{i=1}^{m}u_i' - \sum\limits_{i=1}^{2m} \dfrac{\log(f_i')}{t} $;
		\If{$ F' - F \leq \alpha \cdot s \cdot (\trsp{\nabla F}\cdot\dif\vec{f}) $}
		\State \textbf{break};
		\EndIf
		\EndIf
		\State $ s \gets \beta \cdot s $;
		\EndWhile
		\State $ \vec{r}\gets\vec{r}' $;
		\State $ \vec{u}\gets\vec{u}' $;
		\State $ \vec{f}\gets\vec{f}' $;
		\EndWhile
		
		\State \Return $ \vec{r} $;
		\EndFunction
	\end{algorithmic}
\end{breakablealgorithm}

\subsubsection{Homotopy Method}

The \textit{Homotopy} method exploits the fact that the objective function in \eqref{QP_lambda} undergoes a homotopy from the $ \ell_2 $ constraint to the $ \ell^1 $ objective in \eqref{QP_lambda} as $ \lambda $ decreases \cite{asif2010dynamic,asif2013fast,asif2014homotopy}. The implementations of the homotopy method have been publicly released\footnote{A MATLAB implementation can be found at \url{https://intra.ece.ucr.edu/~sasif/homotopy/index.html}.}. We present the \Algr \ref{alg-homotopy} for solving the  problem \eqref{QP_lambda}  by the the Homotopy method \cite{asif2014homotopy} with the procedure $ \ProcName{L1OptHomotopy}$. For more details about the mathematical principles, please see Appendix \ref{sec:HPL1Opt}.

\begin{breakablealgorithm}
    \caption{Solving the $\ell^1$-norm REV via Homotopy method.}\label{alg-homotopy}
    \begin{algorithmic}[1]
        \Require $ \mat{D}\in\ES{R}{(m-n)}{m},\vec{w}\in\ES{R}{(m-n)}{1},\lambda>0 $, $ \cpvar{maxiter}\in\mathbb{N}^{+} $.
        \Ensure $ \vec{r}\in\ES{R}{m}{1} $.
        \Function{L1OptHomotopy}{$ \mat{D}, \vec{w},\lambda $, $ \cpvar{maxiter} $}
        \State $ \vec{r} \gets \mathbf{0}_{m\times1} $;\Comment{initialize solution}
        \State $ \vec{z} \gets \mathbf{0}_{m\times1} $;
        \State $ \vec{p} \gets -\trsp{\mat{D}} \vec{w}\in \ES{R}{m}{1} $;
        \State $\hat{\vec{p}}\gets \left[\abs{p_1}, \abs{p_2}, \cdots, \abs{p_m}\right]$;
        \State $ \mpair{\scrd{p}{max}}{\vec{\xi}} \gets \ProcName{SearchMax}(\hat{\vec{p}}) $;\Comment{for the maximum $\scrd{p}{max}$ and index vector $ \vec{\xi}=[\xi_1, \cdots, \xi_q]$}
        \State $ \vec{z}[\vec{\xi}] \gets -\text{sign}(\vec{p}[\vec{\xi}]) $;\Comment{primal sign}
        \State $ \vec{p}[\vec{\xi}] \gets \scrd{p}{max} \cdot \text{sign}(\vec{p}[\vec{\xi}]) $;\Comment{dual sign}
        
        \State $ \mat{B} \gets \trsp{[\mat{D}(:, \vec{\xi})]} \mat{D}(:, \vec{\xi}) $;
        
        \While{$ \cpvar{TRUE} $}
        \State \blue{$ \vec{r}_{k} \gets \vec{r} $;}
        \State $ \vec{\gamma} \gets \vec{\xi} $;\Comment{primal support}
        
        \State $ \vec{v} \gets \mathbf{0}_{m\times 1} $;
        \State $ \vec{v}[\vec{\gamma}] \gets \inv{\mat{B}} \cdot \vec{z}[\vec{\gamma}] $;\Comment{update direction}
        \State $ \dif\vec{k} \gets \trsp{\mat{D}} \mat{D} \vec{v} $;
        
        \State $ \seq{\delta, i_{\delta}, i_{\cpvar{shk}}, \cpvar{flag}} \gets \ProcName{CalcStepHomotopy}$($\vec{\gamma}$,  $\vec{\gamma}$,  $\vec{r}_{k}$, $\vec{v}$, $\vec{p}$, $\dif\vec{k}$, $\scrd{p}{max} $);\Comment{\blue{compute the step size}}
        \State $ \vec{r} \gets \vec{r}_{k} + \delta \cdot \vec{v} $;\Comment{update the solution}
        \State $ \vec{p} \gets \vec{p} + \delta \cdot \dif\vec{k} $;
        
        \If{$ (\scrd{p}{max}-\delta) \leq \lambda $}
        \State $ \vec{r} \gets \vec{r}_{k} + (\scrd{p}{max} - \lambda) \cdot \vec{v} $;\Comment{final solution}
        \State \textbf{break}
        \EndIf
        \State $ \scrd{p}{max} \gets \scrd{p}{max} - \delta $;\Comment{update the homotopy parameter}
        
        \If{$ \cpvar{flag} == 1 $}\Comment{remove an element from the support set $\vec{\gamma}$}
        \State $ L \gets \ProcName{Length}(\vec{\gamma}) $;
        \State $ \cpvar{idx} \gets \ProcName{FindIndex}(\vec{\gamma} == i_{\cpvar{shk}}) $;
        \State $ {\gamma}_{\cpvar{idx}} \gets {\gamma}_L $; \Comment{Swap the elements at positions $ \cpvar{idx} $ and $ L $}
        \State $ {\gamma}_L \gets i_{\cpvar{shk}} $;
        \State $ \vec{\xi} \gets \vec{\gamma}(1\!:\!L - 1) $;
        
        \State $ \mat{B} \gets \ProcName{SwapRow}(\mat{B}, \cpvar{idx}, L) $;
        \State $ \mat{B} \gets \ProcName{SwapCol}(\mat{B}, \cpvar{idx}, L) $;
        
        \State $ \mat{B} \gets \mat{B}(1\!:\!L - 1, 1\!:\!L - 1) $;
        \State $ \vec{r}[{i_{\cpvar{shk}}}] \gets 0 $;
        \Else\Comment{add a new element to the support}
        \State $ \vec{\xi} \gets \left[ \begin{matrix}\vec{\gamma}\\ i_{\delta}\end{matrix} \right] $;
        \State $ \mat{C} \gets \trsp{[\mat{D}(:, \vec{\gamma})]} \mat{D}(:, i_{\delta}) $;
        \State $ \mat{B} \gets \left[ \begin{matrix} \mat{B} & \mat{C} \\ \trsp{\mat{C}} & \trsp{[\mat{D}(:, i_{\delta})]} \mat{D}(:, i_{\delta}) \end{matrix} \right] $;
        \State $ \vec{r}[i_{\delta}] \gets 0 $;
        \State $ \vec{\gamma} \gets \vec{\xi} $;
        \EndIf
        
        \State $ \vec{z} \gets \mathbf{0}_{m\times 1} $;\Comment{update primal and dual sign}
        \State $ \vec{z}[\vec{\xi}] \gets -\text{sign}(\vec{p}[\vec{\xi}]) $;
        \State $ \vec{p}[\vec{\gamma}] \gets \scrd{p}{max}\cdot\text{sign}(\vec{p}[\vec{\gamma}]) $;        
        \EndWhile        
        \State\Return $ \vec{r} $;
        \EndFunction
    \end{algorithmic}
\end{breakablealgorithm}

The procedure \ProcName{CalcStepHomotopy} arising in the Line 15 of \Algr \ref{alg-homotopy} is given in
\Algr \ref{alg-homotopy_update_primal}, which is used to calculate the smallest step-size that causes changes in the support set of $ \vec{r} $ or $ \lambda $.

\begin{breakablealgorithm}
    \caption{Calculate the smallest step-size that causes changes in the support set of $ \vec{r} $ or $ \lambda $}.\label{alg-homotopy_update_primal}
    \begin{algorithmic}[1]
        \Require $\vec{\gamma}_{r}$ for the support of $ \vec{r} $,  $\vec{\gamma}_{\lambda}$ for the support of $ \lambda $,  current solution $\vec{r}_{k}$, primal updating direction $\vec{v}$, dual sign vector $\vec{p}$, dual updating direction $\dif\vec{k}$, $\varepsilon>0 $.
        \Ensure Step $ \delta$, index $i_{\delta}$ to be added from the support $ \vec{\gamma}_{\lambda} $, index $i_{\cpvar{shk}}$ to be removed from the support $ \vec{\gamma}_{r} $, $ \cpvar{flag}\in \set{0,1}$ for the  added/removed index.
        \Function{CalcStepHomotopy}{$\vec{\gamma}_{r}$,  $\vec{\gamma}_{\lambda}$,  $\vec{r}_{k}$, $\vec{v}$, $\vec{p}$, $\dif\vec{k}$, $\varepsilon $}
        \State $ \vec{\gamma}_{c}\gets\set{1,2,\cdots,m}-\set{\vec{\gamma}_{\lambda}} $;
        \State $ \vec{\gamma}_{c}\gets\ProcName{Sort}(\vec{\gamma}_{c}) $;\Comment{sorting \& vectorization}
        \State $\vec{\delta}\gets (\epsilon\vec{1} - \vec{p}[\vec{\gamma}_{c}])\oslash(\vec{1} + \dif\vec{k}[\vec{\gamma}_{c}])$;
        \State $\cpvar{idx}_{1}\gets \ProcName{FindIndex}(\vec{\delta}>0)$;\Comment{positive components}
        \If{$\cpvar{idx}_1==\emptyset$}
            \State $\delta_1 \gets +\infty$;
        \Else
            \State $\seq{\delta_1, i_{\delta_{1}}} \gets \ProcName{SearchMin}(\vec{\delta}[\cpvar{idx}_1])$;
        \EndIf
        
        \State $\vec{\delta}\gets(\epsilon\vec{1} + \vec{p}[\vec{\gamma}_{c}])\oslash(\vec{1} - \dif\vec{k}[\vec{\gamma}_{c}])$;
        \State $\cpvar{idx}_2\gets \ProcName{FindIndex}(\vec{\delta}>0)$;
        \If{$\cpvar{idx}_2==\emptyset$}
            \State $\delta_2 \gets +\infty$;
        \Else
            \State $\seq{\delta_2, i_{\delta_{2}}} \gets \ProcName{SearchMin}(\vec{\delta}[\cpvar{idx}_2])$;
        \EndIf
        
        \If{$\delta_1 > \delta_2$}
        \State $\delta \gets \delta_2$;
        \State $i_\delta \gets \vec{\gamma}_{c}[\cpvar{idx}_2[i_{\delta_2}]]$;
        \Else
        \State $\delta \gets \delta_1$;
        \State $i_\delta \gets \vec{\gamma}_{c}[\cpvar{idx}_1[i_{\delta_1}]]$;
        \EndIf
        
        \State $\vec{\delta} \gets -\vec{r}_{k}[\vec{\gamma}_{{r}}]\oslash\vec{v}[\vec{\gamma}_{{r}}]$;
        \State $\cpvar{idx}_3 \gets  \ProcName{FindIndex}(\vec{\delta}>0)$;
        \State $\seq{\delta_3, i_{\delta_3}} \gets \ProcName{SearchMin}(\vec{\delta}[\cpvar{idx}_3])$;
        
        \If{$\delta_3 \leq \delta$}\Comment{an element is removed from support of $\vec{r}$}
        \State $\cpvar{flag} \gets 1$;
        \State $\delta \gets \delta_3$;
        \State $i_{\cpvar{shk}} \gets \vec{\gamma}_{r}[\cpvar{idx}_3[i_{\delta_3}]]$;
        \Else\Comment{a new element enters the support of $\lambda$}
        \State $\cpvar{flag} \gets 0$;
        \State $i_{\cpvar{shk}} \gets -1$;
        \EndIf
        
        \State\Return $ \seq{\delta, i_{\delta}, i_{\cpvar{shk}}, \cpvar{flag}} $;
        \EndFunction
    \end{algorithmic}
\end{breakablealgorithm}

\subsubsection{Iterative Shrinkage-Thresholding Method}

The \textit{Iterative Shrinkage-Thresholding} (IST) method can obtain the REV by solving problem \eqref{QP_lambda}. The implementation of the IST method have been released online \footnote{A MATLAB implementation called \textit{Sparse Reconstruction} by \textit{Separable Approximation} (SpaRSA) is available at \url{http://www.lx.it.pt/~mtf/SpaRSA/}.}. We present the \Algr \ref{alg-IST}  for  solving the problem \eqref{QP_lambda} with the procedure $ \ProcName{L1OptIST}$. The detailed mathematical principles are presented in Appendix \ref{sec:ISTL1Opt}.

\begin{breakablealgorithm}
    \caption{Solving the $\ell^1$-norm REV via Iterative Shrinkage-Thresholding method.}\label{alg-IST}
    \begin{algorithmic}[1]
        \Require $ \mat{D}\in\ES{R}{(m-n)}{m},\vec{w}\in\ES{R}{(m-n)}{1},\lambda>0,\varepsilon>0 $, $ \cpvar{maxiter}\in\mathbb{N}^{+} $.
        \Ensure $ \vec{r}\in\ES{R}{m}{1} $.
        \Function{L1OptIST}{$ \mat{D}, \vec{w},\lambda,\varepsilon $, $ \cpvar{maxiter} $}
        \State $\alpha_{\min} \gets 10^{-30}$;\Comment{set parameters}
        \State $\alpha_{\max} \gets 10^{30}$;
        \State $\vec{r} \gets \vec{0}_{m\times 1}$;\Comment{initialize}
        \State $\vec{s} \gets \mat{D}\vec{r} - \vec{w}$;\Comment{compute residual}
        \State $\alpha \gets 1$;\Comment{initialize}
        \State $f \gets 0.5\cdot\trsp{\vec{s}} \vec{s} + \lambda \cdot \sum\limits_{i=1}^{m}\abs{r_i}$;
        \State $ \cpvar{iter}\gets0 $;
        \While{$\cpvar{iter}<\cpvar{maxiter}$}
        \State $ \cpvar{iter}\gets\cpvar{iter}+1 $;
        \State $\nabla\vec{q} \gets \trsp{\mat{D}}(\vec{s})$;\Comment{compute gradient}
        \State $\vec{r}_{\cpvar{prev}} \gets \vec{r}$;\Comment{save previous values}
        \State $f_{\cpvar{prev}} \gets f$;
        \State $\vec{s}_{\cpvar{prev}} \gets \vec{s}$;
        
        \State $\vec{r} \gets \text{soft}\left(\vec{r}_{\cpvar{prev}} - \dfrac{\nabla\vec{q}}{\alpha}, \dfrac{\lambda}{\alpha}\right)$;
        \State $\dif\vec{r} \gets \vec{r} - \vec{r}_{\cpvar{prev}}$;\Comment{update differences}
        \State $ \vec{z}\gets\mat{D}\cdot\dif\vec{r} $;
        \State $\vec{s} \gets \vec{s}_{\cpvar{prev}} + \vec{z}$;\Comment{update residual}
        \State $f \gets 0.5 \cdot \trsp{\vec{s}} \vec{s} + \lambda \cdot \sum\limits_{i=1}^{m}\abs{r_i}$;
        
        \State $\delta \gets \trsp{(\dif\vec{r})} \cdot \dif\vec{r}$;\Comment{update $\alpha$}
        \State $\gamma \gets \trsp{\vec{z}} \cdot \vec{z}$;
        \State $\alpha \gets \min\left(\alpha_{\max}, \max\left(\alpha_{\min}, \dfrac{\gamma}{\delta}\right)\right)$;
        
        \If{$\dfrac{\abs{f-f_{\cpvar{prev}}}}{f_{\cpvar{prev}}} \leq \varepsilon$}
            \State \textbf{break}
        \EndIf
        \EndWhile
        
        \State\Return $ \vec{r} $;
        \EndFunction
    \end{algorithmic}
\end{breakablealgorithm}

\subsubsection{Alternating Directions Method}

The \textit{Alternating Directions Method} (ADM) can obtain the REV by solving problem \eqref{BPDN}. The implementation of the ADM is also available online \footnote{A MATLAB toolbox of the ADM algorithm named with YALL1 is provided at \url{https://yall1.blogs.rice.edu/}.}. The procedure
\ProcName{L1OptADM} listed in \Algr \ref{alg-ADM} is used to solve the problem \eqref{BPDN} with the Alternating Directions Method \cite{yang2011alternating}. The detailed mathematical principles are presented in Appendix \ref{sec:ADML1Opt}.

\begin{breakablealgorithm}
    \caption{Solving the $\ell^1$-norm REV via Alternating Directions method.}\label{alg-ADM}
    \begin{algorithmic}[1]
        \Require $ \mat{D}\in\ES{R}{(m-n)}{m},\vec{w}\in\ES{R}{(m-n)}{1},\epsilon>0 $, $ \cpvar{maxiter}\in\mathbb{N}^{+} $.
        \Ensure $ \vec{r}\in\ES{R}{m}{1} $.
        \Function{L1OptADM}{$\mat{D}, \vec{w}, \epsilon$, $ \cpvar{maxiter} $}
        
        \State $\zeta \gets 1.618$; \Comment{ADM parameter}
        \State $\mu \gets \dfrac{1}{m-n}\sum\limits_{i=1}^{m-n}\abs{w_i}$;

        \State $\vec{r} \gets \trsp{\vec{D}}\vec{w}$; \Comment{initialization}
        \State $\vec{z} \gets \vec{0}_{m\times1}$;
        \State $\vec{y} \gets \vec{0}_{(m-n)\times1}$;
        \State $\vec{g} \gets \vec{0}_{m\times1}$;
        
        \State $\vec{s} \gets \mat{D}\cdot\left(\vec{g} - \vec{z} + \dfrac{\vec{r}}{\mu}\right) - \dfrac{\vec{w}}{\mu}$;\Comment{calculate step}
        \State $\alpha \gets \dfrac{\trsp{\vec{s}}\vec{s}}{\trsp{\vec{s}}\mat{D}\trsp{\mat{D}}\vec{s}}$;
        
        \State $ \cpvar{iter}\gets0 $;
        \While{$\cpvar{iter}<\cpvar{maxiter}$}
        \State $ \cpvar{iter}\gets\cpvar{iter}+1 $;
        \State $\vec{s} \gets \mat{D}\cdot\left(\vec{g} - \vec{z} + \dfrac{\vec{r}}{\mu}\right) - \dfrac{\vec{w}}{\mu}$;
        
        \State $\vec{y} \gets \vec{y} - \alpha \cdot \vec{s}$;
        \State $\vec{g} \gets \trsp{\vec{D}}\vec{y}$;
        
        \State $\vec{z} \gets \vec{g} + \dfrac{\vec{r}}{\mu}$;
        \For{$ i\in\set{1,2,\cdots,m} $}
            \If{$ \abs{z}_i>1 $}
            \State $ z_i\gets\text{sign}(z_i) $;
            \EndIf
        \EndFor
        \State $\vec{r} \gets \vec{r} + \zeta \cdot \mu \cdot (\vec{g} - \vec{z})$;

        \If{$\cfrac{\normp{\mat{D}\vec{r} - \vec{w}}{2}}{\normp{\vec{w}}{2}} \le \epsilon $}
        \State \textbf{break}
        \EndIf
        \EndWhile
        \State\Return $ \vec{r} $;
        \EndFunction
    \end{algorithmic}
\end{breakablealgorithm}

\subsubsection{Proximity Operator-Based Method}

By constructing an indicator function, the constrained optimization problems can be transformed into a unified unconstrained problem \cite{chen2015fast}. The indicator function of a closed convex set $ \Omega_\epsilon $ is defined as
\begin{equation}\label{POB-1}
	\indcc{\Omega_\epsilon}(\vec{r}) = 
	\begin{cases}
       0,&\vec{r}\in \Omega_\epsilon;\\
		+\infty,&\vec{r}\notin \Omega_\epsilon.
	\end{cases}
\end{equation}
Then the problem \eqref{BPDN} can be expressed by
\begin{equation}\label{POB-2}
	 \min_{\vec{r}\in \ES{R}{m}{1}}\norm{\vec{r}}_1+\indcc{\Omega_\epsilon}(\mat{D}\vec{r}-\vec{w}).
\end{equation}
It is trivial that the problem \eqref{POB-2} reduces to the problem \eqref{BP} when $ \epsilon=0 $. The proximity operator $ \Prox_{\indcc{\Omega_\epsilon} } $ for \eqref{POB-1} is the \textit{soft-thresholding} operator defined by
\begin{equation*}\label{POB-4}
	\begin{aligned}
		\Prox_{\indcc{\Omega_\epsilon}}(\vec{y},\vec{w})=\vec{w}+\min\set{\!1,\frac{\epsilon\!}{\normp{\vec{y}-\vec{w}}{2}}}\cdot(\vec{y}-\vec{w}).
	\end{aligned}
\end{equation*}
There exists the following iterative scheme
\begin{equation}\label{POB-5}
	\left\{\begin{aligned}
		&\begin{aligned}
			\vec{r}^{(k+1)}=\textrm{soft}&\left(\left(\mat{I}_m-\frac{\tau}{\mu}\trsp{\mat{D}}\mat{D}\right)\vec{r}^{(k)}-\right.\\
			&\left.\frac{\tau}{\mu}\trsp{\mat{D}}\left(\vec{y}^{(k)}-\vec{u}^{(k)}\right),\frac{1}{\mu}\right)
		\end{aligned},\\
		&\vec{u}^{(k+1)}=\Prox_{\indcc{\Omega_\epsilon}}\left(\mat{D}\vec{r}^{(k+1)}+\vec{y}^{(k)},\vec{w}\right),\\
		&\vec{y}^{(k+1)}=\mat{D}\vec{r}^{(k+1)}+\vec{y}^{(k)}-\vec{u}^{(k+1)}.
	\end{aligned}\right.
\end{equation}
where $ \tau>\mu\normp{\mat{D}}{2}^2>0 $ and the soft-thresholding function is defined in \eqref{IST-5}. \Algr \ref{alg-BPiterative} is a simplified implementation of \eqref{POB-5} by eliminating the intermediate variable $ \vec{u}^{(k)} $ \cite{chen2015fast}.
\begin{breakablealgorithm}
	\caption{Solving the $\ell^1$-norm REV via POB method.}\label{alg-BPiterative}
	\begin{algorithmic}[1]
		\Require $ \mat{D}\in\ES{R}{(m-n)}{m},\vec{w}\in\ES{R}{(m-n)}{1},\epsilon>0,\tau>0,\mu>0 $ and $ \tau>\mu\normp{\mat{D}}{2}^2 $, $ \cpvar{maxiter}\in\mathbb{N}^{+} $.
		\Ensure $ \vec{r}\in\ES{R}{m}{1} $.
		\Function{L1OptPOB}{$ \mat{D}, \vec{w},\epsilon,\tau,\mu $, $ \cpvar{maxiter} $}
		\State $ \vec{y}\gets\vec{0}_{m-n} $;
		\State $ \vec{r}\gets\vec{0}_{m} $;
		\State $ \vec{z}\gets\vec{y}-(\mat{D}\vec{r}-\vec{w}) $;
		\State $ \cpvar{iter}\gets0 $;
        \While{$\cpvar{iter}<\cpvar{maxiter}$}
        \State $ \cpvar{iter}\gets\cpvar{iter}+1 $;
			\State $ \vec{s}\gets\vec{r} $;
			\State $ \vec{t}\gets\vec{s}-\dfrac{\mu}{\tau}\trsp{\mat{D}}(2\vec{y}-\vec{z}) $;
			\For{$ i\in\set{1,\cdots,m} $}\Comment{soft-thresholding for $r_i$}
				\If{$ \abs{t_i}>\dfrac{1}{\tau} $}
					\State $ r_i\gets \left(\abs{t_i}-\dfrac{1}{\tau}\right)\cdot\sign(t_i) $;
				\Else
					\State $ r_i\gets 0 $;
				\EndIf
			\EndFor
			\If{$ \dfrac{\normp{\vec{r}-\vec{s}}{2}}{\normp{\vec{s}}{2}}<10^{-6} $}
				\State \textbf{break}
			\EndIf
			\State $ \vec{z}\gets\vec{y} $;
			\State $ \vec{t}\gets\mat{D}\vec{r}+\vec{z}-\vec{w} $;
			\If{$ \normp{\vec{t}}{2}\le\epsilon $}
				\State $ \vec{y}\gets\vec{0}_{m-n} $;
			\Else
				\State $ \vec{y}\gets\left(1-\dfrac{\epsilon}{\normp{\vec{t}}{2}}\right)\vec{t} $;
			\EndIf
		\EndWhile
		\State\Return $ \vec{r} $;
		\EndFunction
	\end{algorithmic}
\end{breakablealgorithm}

\subsection{Unified Framework  for $ \ell^1 $-norm Approximation via Minimizing the $ \ell^1 $-norm of REV}

According to the \Theo \ref{fzq-theorem-1}, we can design a unified framework for $ \ell^1 $-norm approximation via minimizing the $ \ell^1 $-norm of REV. We now give the pseudo-code for the $ \ell^1 $-norm approximation via minimizing $ \ell^1 $-norm residual vector, please see \Algr \ref{alg-L1residual}.
\begin{breakablealgorithm}
	\caption{Solving the $ \ell^1 $-norm approximation to the multivariate linear model by minimizing the $ \ell^1 $-norm of REV.}\label{alg-L1residual}
	\begin{algorithmic}[1]
		\Require the matrix $ \mat{A}\in\ES{R}{m}{n}$ and vector $\vec{b}\in\ES{R}{m}{1}$ such that $m>n\ge2 $, positive number $\epsilon>0$, function object $\cpvar{f}$ with variable list of arguments specified by the grammar
\begin{align*}
\cpvar{f}: \ES{R}{(m-n)}{m}\times \ES{R}{(m-n)}{1}\times \mathbb{R}^+\times \cdots
&\to \ES{R}{m}{1} \\
(\mat{D},\vec{w},\epsilon,...) &\mapsto \vec{r}
\end{align*}		
		 for the procedures \ProcName{L1OptLinProg}, ..., 
	     \ProcName{L1OptADM} with the input argument list $(\mat{D},\vec{w},\epsilon)$ and the procedure 
\ProcName{L1OptPOB} with default values $\tau = 0.02$ and $\mu = 0.999\tau/\normp{\mat{D}}{2}^2$ for the extra arguments $\tau$ and $\mu$.
		\Ensure $ \vec{x}\in\ES{R}{n}{1} $.
		\Function{L1ApproxViaMinREV}{$ \cpvar{f},\mat{A}, \vec{b},\epsilon,...$}
		\State $[m,n]\gets \ProcName{Size}(\mat{A})$;\Comment{set the size of  $ \mat{A}\in \ES{R}{m}{n}$}
		\State $ \mat{A}_1\gets\mat{A}(1:n,1:n) $;
		\State $ \mat{A}_2\gets\mat{A}(n+1:m,1:n) $;
		\State $ \mat{A}_1^\dagger\gets\ProcName{ClacInvMatMP}(\mat{A}_1) $;
		\State $ \mat{D}\gets\mybmatrix{-\mat{A}_2\ginv{\mat{A}_1}, \mat{I}_{m-n}} $;
		\State $ \vec{w}\gets\mat{A}_2\mat{A}_1^\dagger\vec{b}(1:n) - \vec{b}(n+1:m) $;
		\State $ \vec{r}\gets\cpvar{f}(\mat{D},\vec{w}, \epsilon,...) $;\Comment{$\ell^1$ optimization}
		\State $ \mat{A}^\dagger\gets\ProcName{ClacInvMatMP}(\mat{A}) $;
		\State $ \vec{x}\gets\mat{A}^\dagger(\vec{b}+\vec{r}) $;\Comment{by \eqref{residual-4}}
		\State\Return $ \vec{x} $;
		\EndFunction
	\end{algorithmic}
\end{breakablealgorithm}

Please note that the \verb|...| in the line 1 and line 8 of \Algr \ref{alg-L1residual} represent the optional arguments $\tau$ and $\mu$, which are obtained from the results of \cite{chen2015fast}. Obviously, the technique of variable list of arguments for designing procedures in the sense of computer programming is taken since the procedure \ProcName{L1OptPOB} needs two extra arguments $\tau$ and $\mu$. The solver \ProcName{ClacInvMatMP} for the Moore-Penrose inverse used in line 5 of \Algr \ref{alg-L1residual} can be implemented using different techniques, such as the Greville column recursive algorithm \cite{GolubLoan2013,ZhangXD2017Matrix}.

\section{Performance Evaluation for the Algorithms}
\label{sec-alg-performance}

In this section, we will compare the performance of the $ \ell^1 $-norm approximation based on the improved minimal $ \ell^1 $-norm residual vector solver with that of linear programming and perturbation methods, in terms of accuracy and running time. It should be noted that our testing platform has the following configuration: Ubuntu 22.04.3 LTS (64-bit); Memory, 32GB RAM; Processor, 12th Gen Intel$ ^\circledR $ Cor$ ^\text{\scriptsize TM} $ i7-12700KF $ \times $ 20; GNU Octave, version 6.4.0.

The precision parameter in the experiment is set to $ 10^{-8} $, and $ \cpvar{maxiter} $ is set to $ 10000 $ for all algorithms except algorithm $ \ProcName{L1ApproxPertCBS} $, which is set to $ 15 $.

\subsection{Sparse Noisy Data and  Redundancy Level}

Let $ n $ be a fixed value and 
\begin{equation}
\DRL = \frac{m}{n}
\end{equation}
be the \textit{data redundancy level} (DRL) of the noisy linear system $ \mat{A}\vec{x}=\vec{b}+\vec{q}$, where $ \vec{q} $ is a sparse noise vector. The \textit{sparsity ratio} of $ \vec{q} $ is defined by
the ratio of the number of non-zero elements to the total length of the sequence, viz.
\begin{equation}
\scrd{\gamma}{sp}(\vec{q}) = \frac{\abs{\set{i: 1\le i\le m, q_i \neq 0 }}}{m}, \quad \vec{q}\in \ES{R}{m}{1}
\end{equation} 


\subsection{Noise-Free/Well-determined Case}\label{experiment-1}

Construct the coefficient matrix $ \mat{A}\in\ES{R}{m}{n} $ and vector $ \vec{p}\in\ES{R}{n}{1} $, where $ m=2^{8},n=2^{7} $ and each element in $ \mat{A} $ and $ \vec{p} $ follows a standard normal distribution. Let $ \vec{b}=\mat{A}\vec{p}\in\ES{R}{m}{1} $, and $ \vec{p} $ is the exact solution to problem \eqref{P_1} \cite{Fzq2022}.

In each experiment, a linear system $ \mat{A}\vec{p}=\vec{b} $ is constructed and the aforementioned algorithm is applied to solve it. This yields the empirical solution $ \scrd{\hat{\vec{p}}}{\Alg}^i $ (where $ \epsilon $ and $ \lambda $ for each algorithm are set to $ 10^{-8} $). The relative error for each experiment is defined as follows
\begin{equation*}\label{Exp-1}
	\eta_\Alg^i=\dfrac{\normp{\scrd{\hat{\vec{p}}}{\Alg}^i-\vec{p}}{2}}{\normp{\vec{p}}{2}}\times100\%.
\end{equation*}
Then, the average relative error is calculated by
\begin{equation*}\label{Exp-2}
	\eta_\Alg=\frac{1}{N}\sum_{i=1}^{N}\eta_\Alg^i,
\end{equation*}
where $ N=30 $ is the number of repetitions. Consequently, for the  algorithm labeled by \Alg, we can compare their relative errors and record an average running time to evaluate their performance. For the convenience of reading, we give some labels for the algorithms discussed, please see \Tab \ref{tab-algr-label}.

\begin{table*}[h]
\caption{Labels for the $\ell^1$-norm optimization method for MLM.} \label{tab-algr-label}
\begin{tabular}{cllll}
\toprule
\textbf{No.} & \textbf{Label} & \textbf{Algorithm} & \textbf{Procedure Name} & \textbf{Argument} $ \cpvar{f} $ \\
\midrule
1& L1-PTB  & \Algr \ref{alg-Perturbation} & \ProcName{L1ApproxPertCBS} & ----- \\
2& L1-LP & \Algr \ref{alg-L1linprog} & \ProcName{L1ApproxLinProg} & -----  \\
3& L1-RES  & \Algr \ref{alg-BPlinprog},\ref{alg-L1residual}  & \ProcName{L1ApproxViaMinREV}(\cpvar{f},...) & \ProcName{L1OptLinProg} \\
4& L1-GPSR  & \Algr \ref{alg-GPSRBB},\ref{alg-L1residual}  & \ProcName{L1ApproxViaMinREV}(\cpvar{f},...) & \ProcName{L1OptGPSR} \\
5& L1-TNIPM & \Algr \ref{alg-TNIPM},\ref{alg-L1residual} & \ProcName{L1ApproxViaMinREV}(\cpvar{f},...) & \ProcName{L1OptTNIPM} \\
6& L1-HP & \Algr \ref{alg-homotopy},\ref{alg-L1residual} & \ProcName{L1ApproxViaMinREV}(\cpvar{f},...) & \ProcName{L1OptHomotopy} \\
7& L1-IST & \Algr \ref{alg-IST},\ref{alg-L1residual} & \ProcName{L1ApproxViaMinREV}(\cpvar{f},...) & \ProcName{L1OptIST} \\
8& L1-ADM & \Algr \ref{alg-ADM},\ref{alg-L1residual} & \ProcName{L1ApproxViaMinREV}(\cpvar{f},...) & \ProcName{L1OptADM} \\
9& L1-POB & \Algr \ref{alg-BPiterative},\ref{alg-L1residual} & \ProcName{L1ApproxViaMinREV}(\cpvar{f},...) & \ProcName{L1OptPOB} \\
\bottomrule
\end{tabular}
\end{table*}
  
\begin{figure}[!htbp]
	\centering
    \subfigure[Relative error for each approximation]{
        	\includegraphics[width=.95\linewidth]{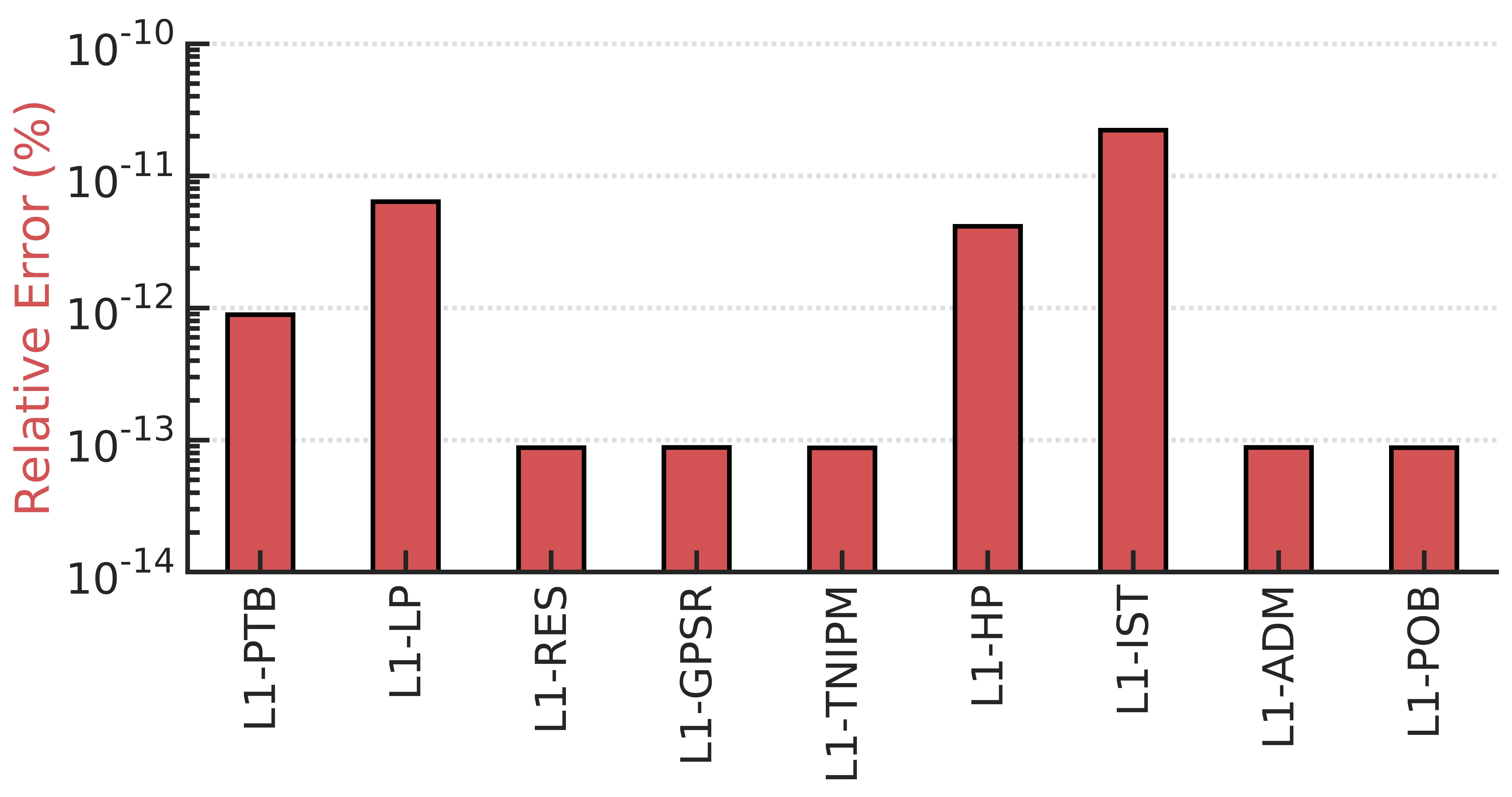}
        	\label{fig-relative-err-noise-free}
     }
    \subfigure[Running time for each approximation]{
        \includegraphics[width=.95\linewidth]{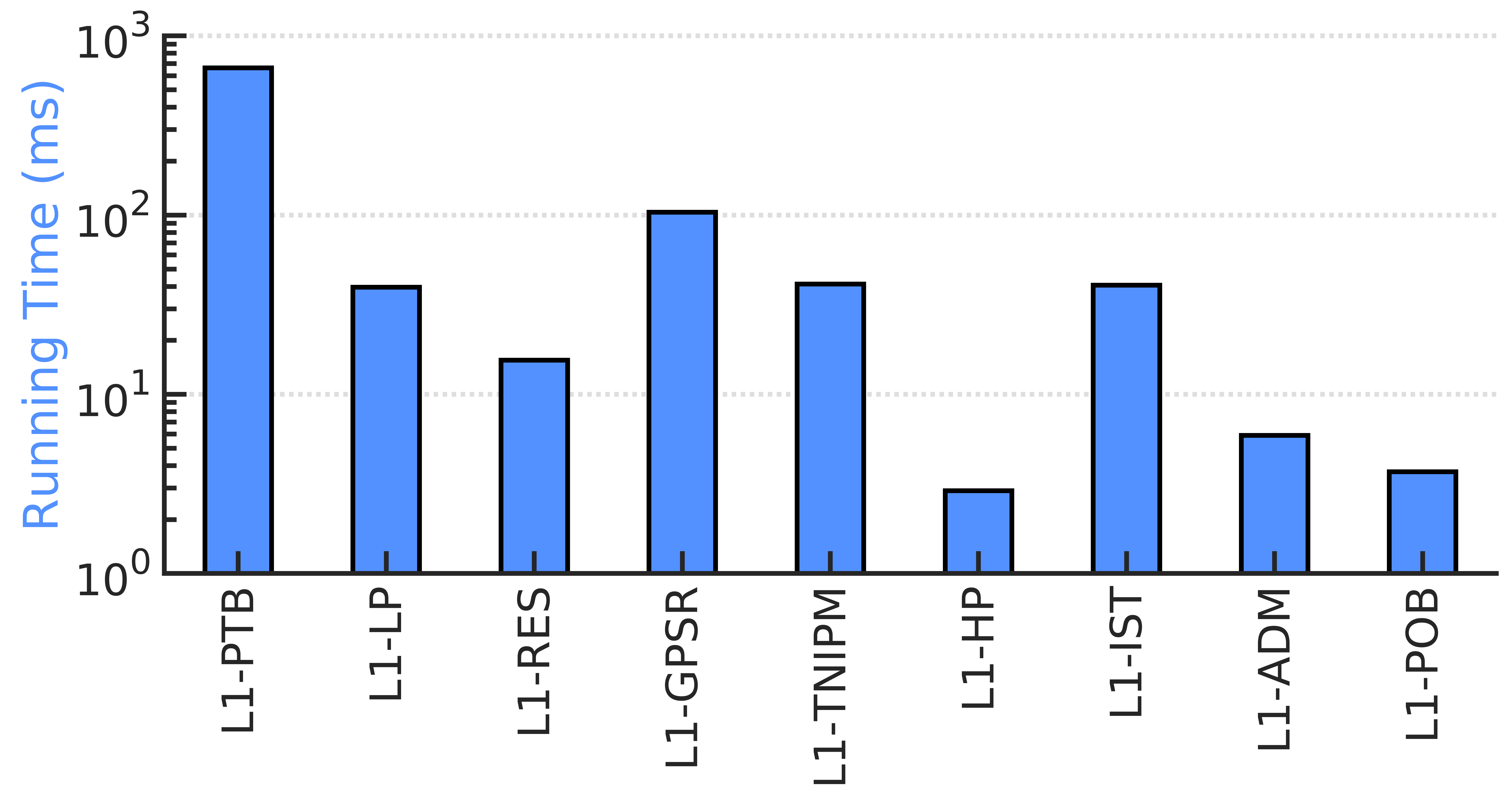}
        \label{fig-running-time-noise-free}
    }
	\caption{Comparison of $\ell^1$-norm optimization method in the noise-free case} \label{noise-free}
\end{figure}

\Fig \ref{noise-free} indicates that all algorithms exhibit high accuracy $ (<10^{-12}) $ and operational efficiency in the noise-free case. All the labels mentioned in the figure can be found in \Tab\ref{tab-algr-label}.

\subsection{Noisy Case}

\subsubsection{Sparse Noise Case}\label{experiment-2}

In each experiment mentioned in section \ref{experiment-1}, the observation vector $ \vec{b} $ is corrupted by the sparse noise $ \vec{q}\in\ES{R}{m}{1} $, thus we have$ \mat{A}\vec{p}=\vec{b}+\vec{q} $. The sparsity ratio of $ \vec{q} $ is specified by
$$ \scrd{\gamma}{sp}(\vec{q}) \in \set{0.25, 0.50, 0.75}$$
with its non-zero elements following a Gaussian distribution with zero mean and variance $ 0.25 $.
\begin{figure}[htb]
	\centering
	\subfigure[Relative error for each approximation]{
		\includegraphics[width=0.95\linewidth]{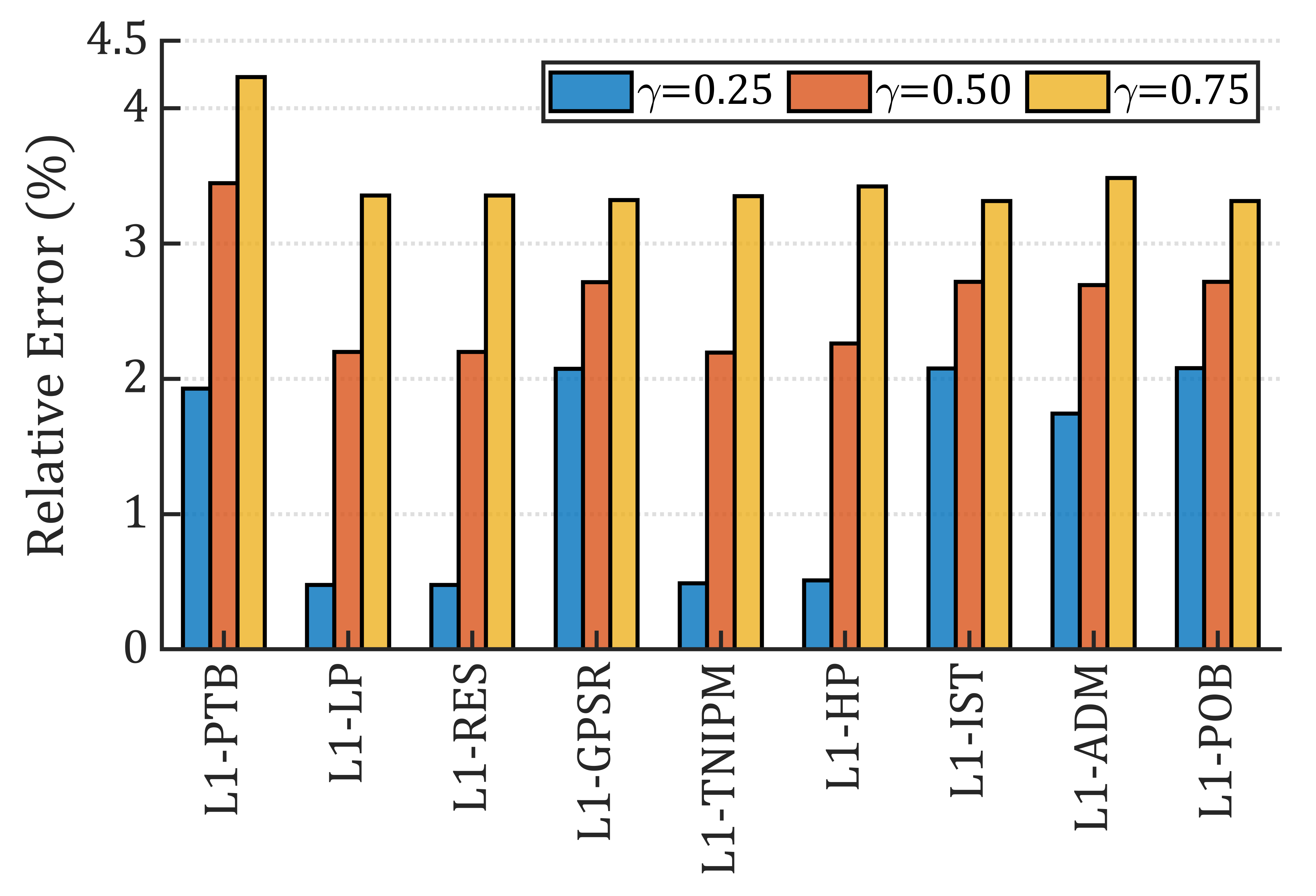}
		\label{fig-relative-err-noise}
	}
    \subfigure[Running time for each approximation]{
        \includegraphics[width=0.95\linewidth]{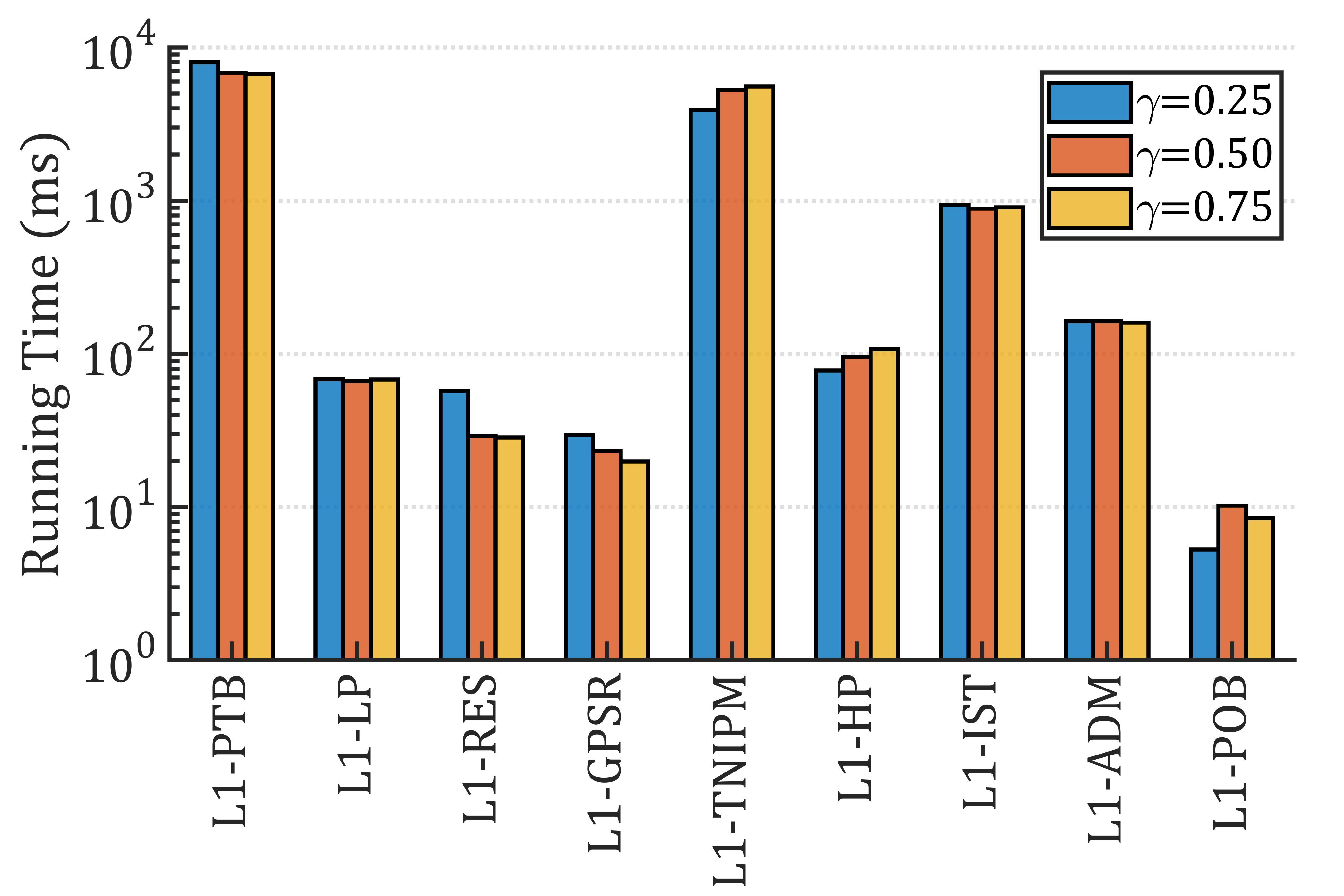}
        \label{fig-running-time-noise}
    }
	\caption{Comparison of $\ell^1$-norm optimization method in the sparse-noise case with sparsity ratio $\scrd{\gamma}{sp}(\vec{q})\in\set{0.25, 0.50, 0.75}$}\label{noise-sparse}
\end{figure}

\Fig \ref{noise-sparse} illustrates the relative error and running time for the algorithms discussed in this paper. From the \Fig \ref{noise-sparse}, we can find some interesting results:
\begin{itemize}
	\item the introduction of noise has an impact on the estimation accuracy and running time of different algorithms;
	\item the impact of noise at different sparsity ratios on the algorithm's running time is relatively weak;
	\item a higher sparsity ratio will increase the relative error of algorithmic estimation;
	\item different solvers have different precision and computational complexity of time.
\end{itemize}

\subsubsection{Impact of Data Redundancy Level}

\begin{figure*}[htb]
	\centering
	\subfigure[Relative error of each algorithm]{
		\includegraphics[width=0.48\textwidth]{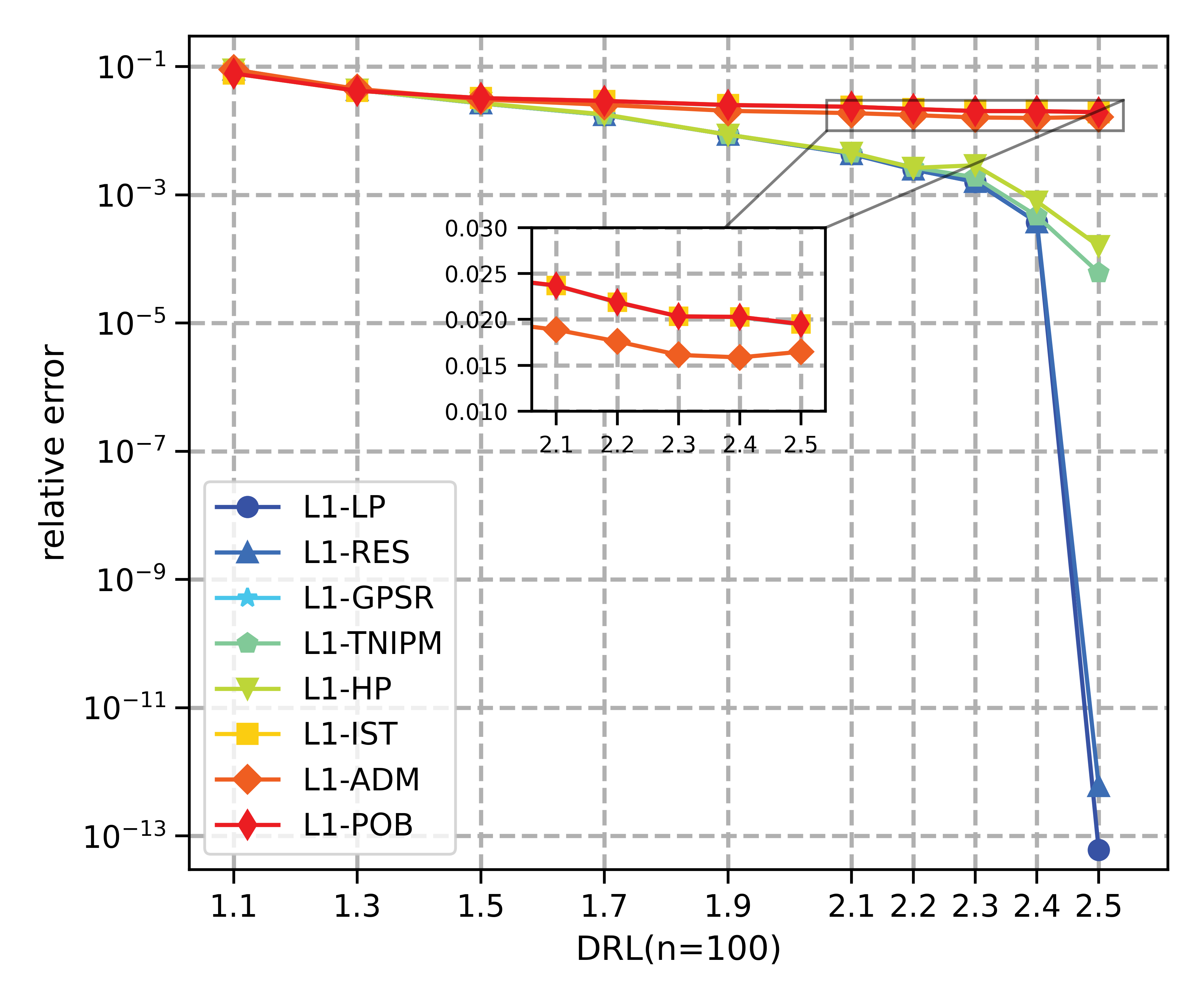}
		\label{fig-DRL-err100}
	}\hfill
	\subfigure[Running time of each algorithm]{
		\includegraphics[width=0.48\textwidth]{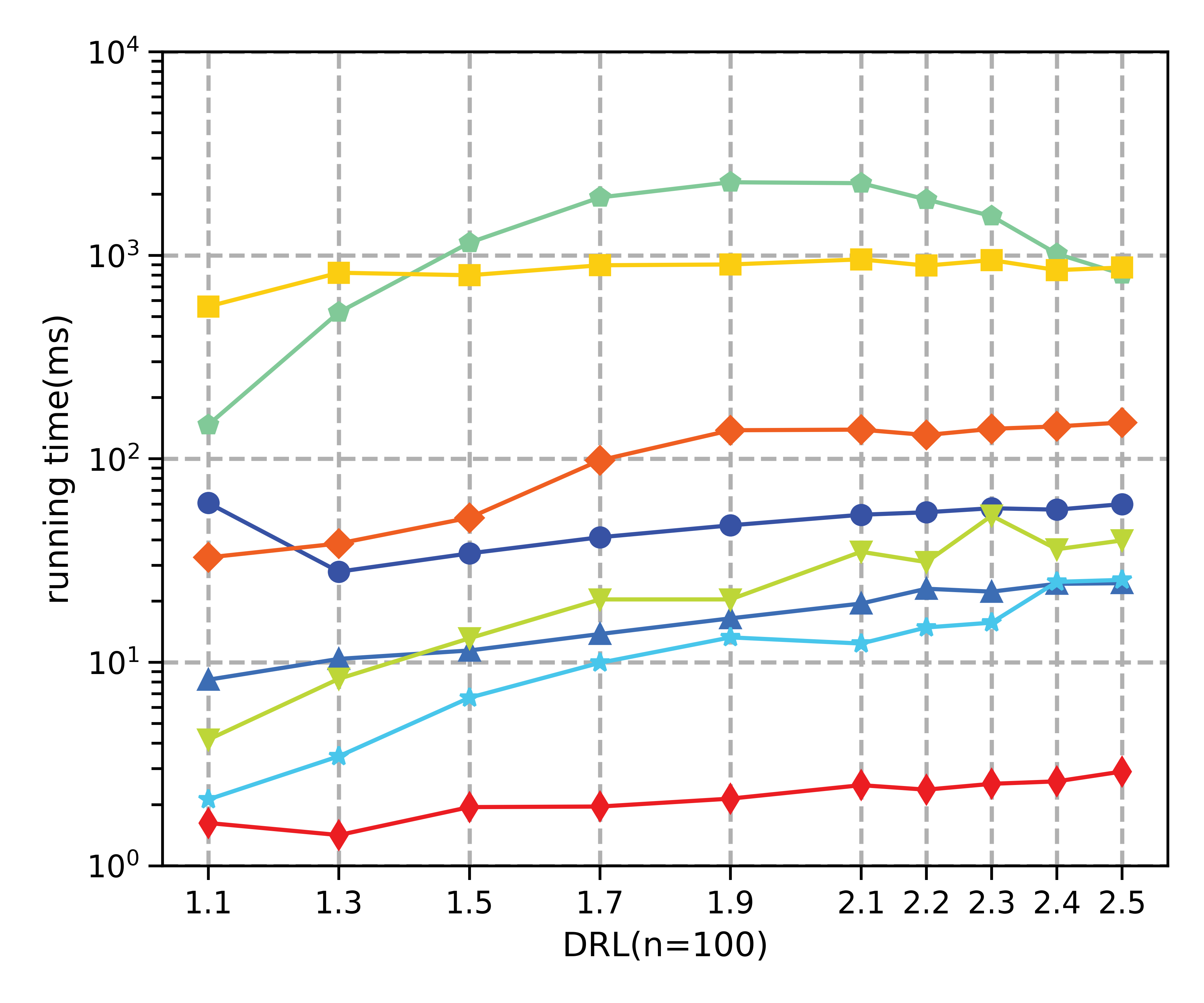}
		\label{fig-DRL-time100}
	}
	\caption{Comparison of each approximation with different data redundancy levels when $ n=100 $}\label{DRL100}
\end{figure*}

\begin{figure*}[htb]
	\centering
	\subfigure[Relative error of each algorithm]{
		\includegraphics[width=0.48\textwidth]{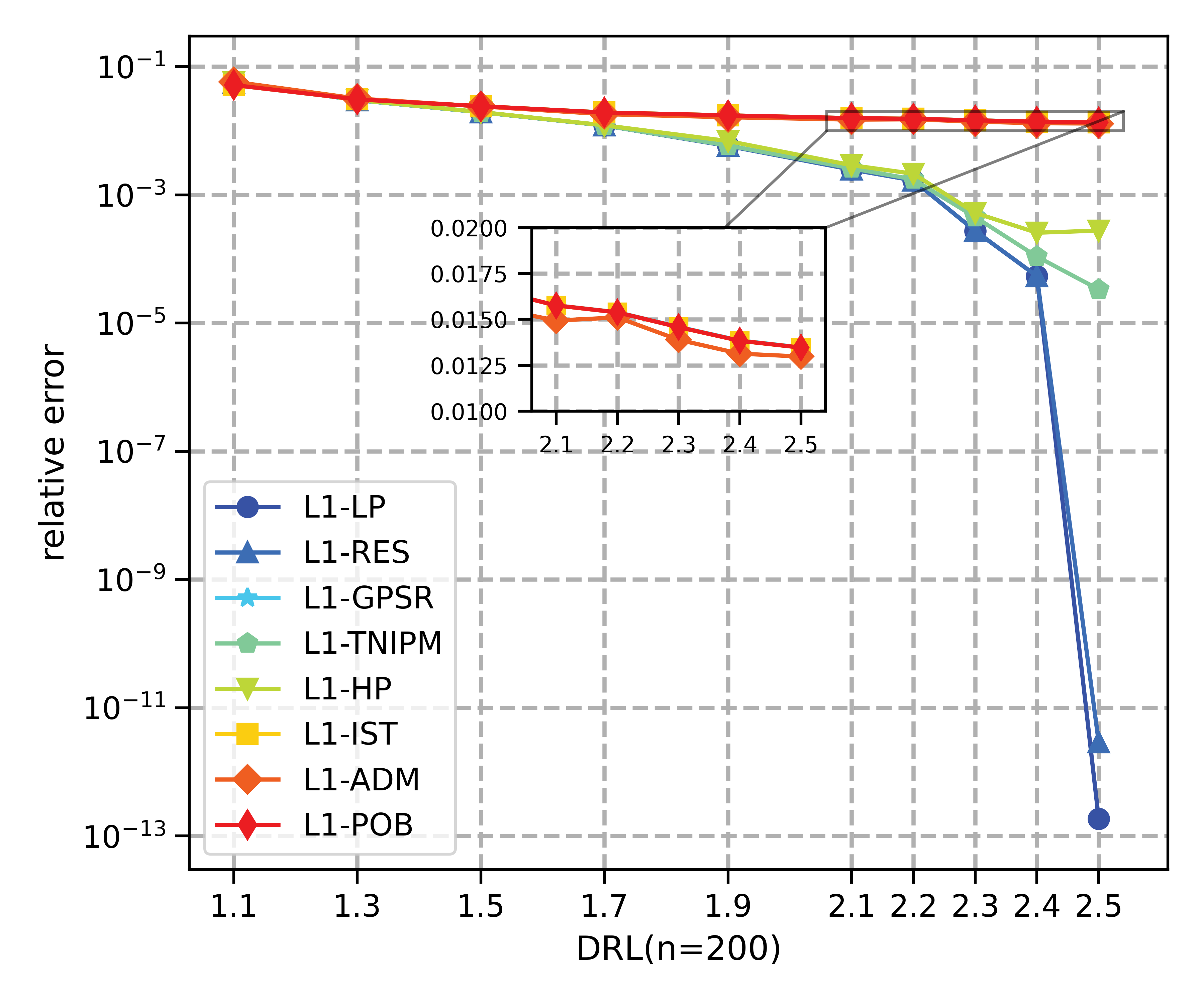}
		\label{fig-DRL-err200}
	}\hfill
	\subfigure[Running time of each algorithm]{
		\includegraphics[width=0.48\textwidth]{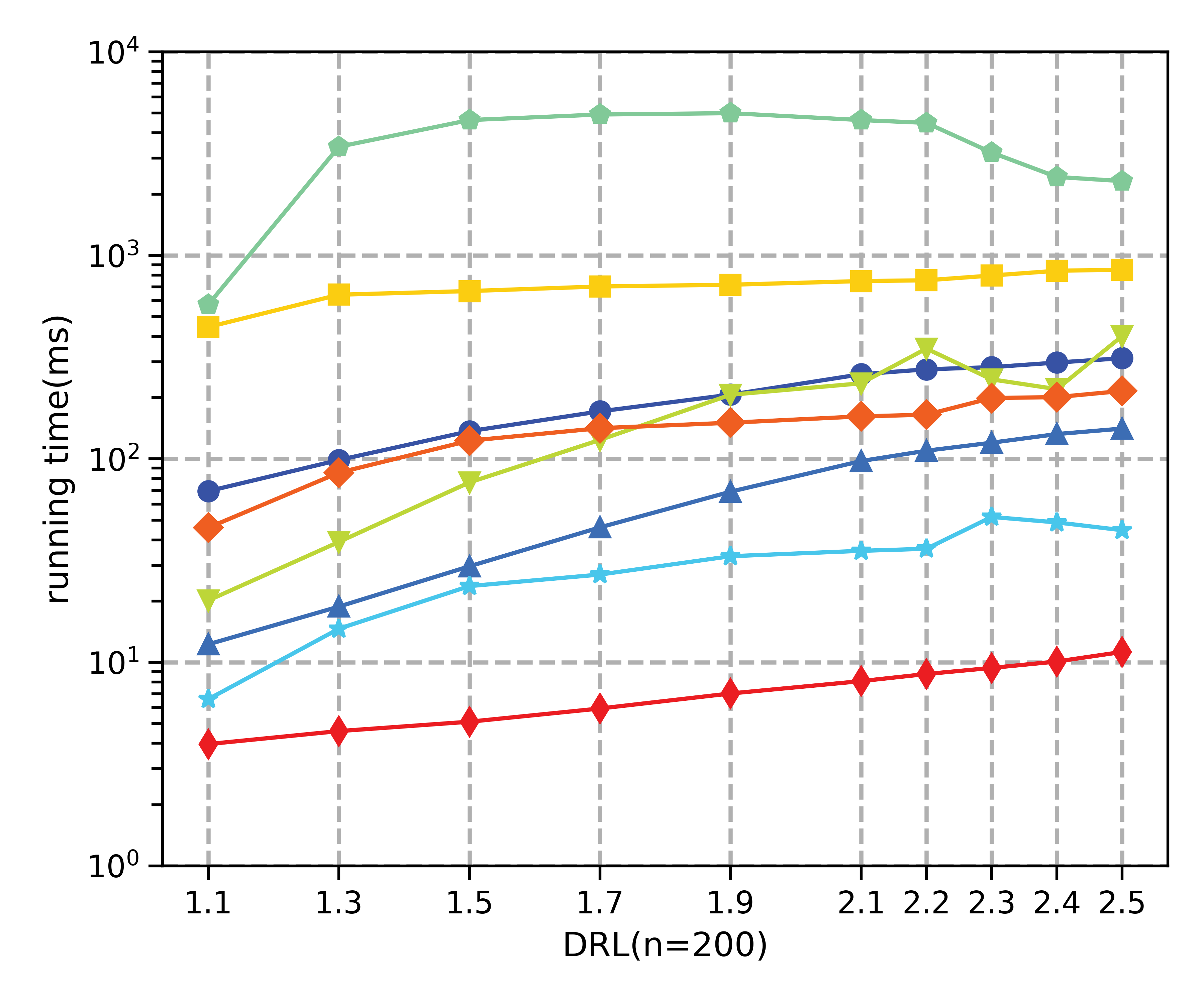}
		\label{fig-DRL-time200}
	}
	\caption{Comparison of each approximation with different data redundancy levels when $ n=200 $}\label{DRL200}
\end{figure*}

\Fig \ref{DRL100} and \Fig \ref{DRL200} depict the correlation between the algorithm performances and the DRL: 
\begin{itemize}
	\item as the DRL increases, the relative errors of different algorithms gradually decrease, and this trend remains unaffected by the total amount of data $ (m\times n) $;
	\item the running time of each algorithm increases as the DRL and the total amount of data grow;
	\item once the DRL reaches a certain value, the accuracy of two linear programming-related algorithms experiences a sudden decreasing.
\end{itemize}

It is important to note that in \Fig \ref{fig-DRL-err100} and \Fig \ref{fig-DRL-err200}, the relative error curves of the L1-GPSR algorithm and the L1-SpaRSA algorithm almost coincide with the L1-POB algorithm. The perturbation approximation method specified by the \Algr \ref{alg-Perturbation} given in the appendices does not perform well in this case and  it is not shown in \Fig \ref{DRL100} and \Fig \ref{DRL200}. 

\section{Conclusions} \label{sec-conclusion}

Solving the problem \eqref{P_1} accurately and efficiently is a challenging task, especially in real-time applications such as the SLAM in computer vision and autonomous driving. The key theoretical contribution of this work is the equivalence theorem for the structure of $\ell_1$ approximation solution to the MLM $\mat{A}\vec{x}=\vec{b}$: the optimal solution $\scrd{\vec{x}}{opt}$ is sum of the traditional LS solution $\mat{A}^\dag \vec{b}$ and the correction term $\mat{A}^\dag\scrd{\vec{r}}{opt}$ specified by the ML1-REV.

The simulations show that our $ \ell^1 $-norm approximation algorithms based on existing $ \ell_1 $-norm optimization algorithms to minimizing the residual vector, {namely the L1-POB, L1-GPRS, L1-HP, ...}, have the following advantages in solving the problem \eqref{P_1}:
\begin{itemize}
	\item[1)] Both the L1-POB and the L1-GPRS have low time complexity. Given the computational platform, the ratio of the running time for the L1-HP algorithm and L1-RES algorithm is approximately $1/5$ if there is no noise. By comparison, the ratio of the running time for the L1-POB algorithm and the L1-RES algorithm is about $1/9$ if the noise is sparse;
	\item[2)] The accuracy of the solutions obtained by our algorithms is comparable to accuracy of the L1-RES algorithm in various scenarios; 
	\item[3)] The implementations are simple and they does not depend on
	any built-in mathematical programming solvers;
 	\item[4)] The L1-POB and L1-GPSR algorithms can not only find the solution with  high accuracy but also with high efficiency when compared with the L1-RES algorithm 
 	in the scenarios with different levels of data redundancy.
\end{itemize}

Among the evaluated algorithms, the L1-POB algorithm distinguishes itself with its low time complexity, which implies that it is suitable for real-time computation in data-constrained scenarios. In the sense of accuracy, all of the $\ell_1$ optimization algorithm for solving the problem \eqref{P_1} exhibit a relative error within $ 3\% $. Particularly, the L1-HP algorithm could offer better accuracy in applications. In the case of high data redundancy, both the L1-LP algorithm and the L1-RES algorithm demonstrate high accuracy exceptionally. By relaxing the constraints, the space of the feasible solutions to the MLM can be expanded, which allows  a broader range of the potential solutions. 

It is still a challenging problem that how to solve the MLM with high accuracy, high efficiency and high robustness via the $ \ell^1 $-norm optimization in computational mathematics. We believe that our explorations and algorithms  could offer new directions in addressing the problem \eqref{P_1} in the future.

For the convenience of usage, all of the algorithms discussed in this paper are implemented with the popular programming languages Python and Octave/Octave,  and the source code has been released on GitHub. We hole that this study could speed up the propagation and adoption of the $\ell_1$-norm optimization in various applications where the MLM appears naturally. 

\subsection*{Acknowledgments} 
This work was supported in part by the National Natural Science Foundation of China under grant numbers 62167003 and 52302421, and in part by the Program of Tianjin Science and Technology Plan under grant number 23JCQNJC00210, and in part by the Research Project on Education and Teaching Reform in Higher Education System of Hainan Province under grant number Hnjg2025ZD-28.

\subsection*{Code Availability}

The code for the implementations of the algorithms dis-
cussed in this paper can be downloaded from the following
GitHub website 
\begin{center}
\url{https://github.com/GrAbsRD/MlmEstAlgorL1} 
\end{center}
For the convenience of easy usage, both Python and Octave/MATLAB codes are provided, please check the following packages:
\begin{itemize}
\item  \verb|MLM-L1Opt-Solver-Matlab-code.zip| 
\item  \verb|MLM-L1Opt-Solver-Python-code.zip| 
\end{itemize} 
Please note that:
\begin{itemize}
\item the Octave/MATLAB code depends on some available solvers mentioned in the footnotes in this paper;
\item the Python code is completely implemented by the authors and it do not depend on any specific solvers available online. 
\end{itemize}

\subsection*{Data Availability}

The data set supporting the results of this study is also available
on the GitHub website 
\begin{center}
\url{https://github.com/GrAbsRD/MlmEstAlgorL1} 
\end{center}
The Octave/MATLAB  script file \verb|GenTestData.m| is used to generate the test data for validating the algorithms discussed in this paper. Three data sample are provided in the \verb|*.txt| files and more can be generated by running the script file \verb|GenTestData.m|. 

\begin{appendix}
\section{Mathematical Principles of Typical $ \ell_1 $-norm Optimization Methods}

\subsection{Linear Programming Method}\label{sec:LinprogL1Opt}

 Farebrother \cite{farebrother2013l1} pointed out that the unconstrained $ \ell^1 $-norm optimization problem $ \min\limits_{\vec{r}\in\ES{R}{m}{1}}\normp{\vec{r}}{1} $ is equivalent to 
\begin{equation}\label{BP-linprog-1}
\min\limits_{\vec{r}^{+},\vec{r}^{-}\in\ES{R}{m}{1}}\trsp{\vec{1}}_m(\vec{r}^{+}+\vec{r}^{-}),\ \st \begin{cases}
\vec{r}^{+}-\vec{r}^{-}=\vec{r}\\
\vec{r}^{+},\vec{r}^{-}\succcurlyeq 0
\end{cases}
\end{equation}
Let
\begin{equation}\label{BP-linprog-2}
\vec{t}=\vec{1}_{2m}=\begin{bmatrix}
\vec{1}_m\\ \vec{1}_m
\end{bmatrix},\quad \vec{\beta}=\begin{bmatrix}
\vec{r}^{+}\\ \vec{r}^{-}
\end{bmatrix}, \quad \vec{\Phi}=\begin{bmatrix}
\vec{D},-\vec{D}
\end{bmatrix}
\end{equation}
and substitute \eqref{BP-linprog-2} into \eqref{BP-linprog-1}, we immediately have the equivalent description \cite{chen2001atomic}
\begin{equation}\label{LP_2}\tag{$ \rm{LP}_2 $}
\min\limits_{\vec{\beta}\in\ES{R}{2m}{1}}\trsp{\vec{t}}\vec{\beta}, \quad  \st\begin{cases}
\vec{\Phi}\vec{\beta}=\vec{w}\\
\vec{\beta}\succcurlyeq \vec{0}
\end{cases}
\end{equation}
where $ \mat{D} $ and $ \vec{w} $ are derived from \eqref{residual-1} and \eqref{residual-2}. Suppose the solution to the  problem \eqref{LP_2} is 
\begin{equation}
\scrd{\vec{\beta}}{opt}=\trsp{[
	\scrd{\vec{r}}{opt}^{+}, \scrd{\vec{r}}{opt}^{-}]}, 
\end{equation}
then the solution to the problem \eqref{BP} must be
\begin{equation}\label{BP-linprog-4}
\scrd{\vec{r}}{opt}=\scrd{\vec{r}}{opt}^{+}-\scrd{\vec{r}}{opt}^{-}=\scrd{\vec{\beta}}{opt}(1:m)-\scrd{\vec{\beta}}{opt}(m+1:2m).
\end{equation}

\subsection{Gradient Projection Method}\label{sec:GPL1Opt}

The first \textit{Gradient Projection} (GP) method is the GPSR method \cite{figueiredo2007gradient}.  By segregating the positive coefficients $ \vec{r}^{+} $ and the negative coefficients $ \vec{r}^{-} $ in $ \vec{r} $, the equation \eqref{QP_lambda}  can be reformulated into the form of standard quadratic programming, namely 
\begin{equation}\label{GPSR-1}
\min\limits_{\vec{z} \in \ES{R}{2m}{1}}Q(\vec{z})=\trsp{\vec{\gamma}}\vec{z}+\frac{1}{2}\trsp{\vec{z}}\mat{B}\vec{z},\ \st\ \vec{z}\succcurlyeq0
\end{equation}
where 
\begin{equation} \label{eq-z-gamma-B}
\vec{z}=\begin{bmatrix}
\vec{r}^{+}\\\vec{r}^{-}
\end{bmatrix},\vec{\gamma}=\begin{bmatrix}
\lambda-\trsp{\mat{D}}\vec{w}\\\lambda+\trsp{\mat{D}}\vec{w}
\end{bmatrix},  \mat{B}=\begin{bmatrix}
\trsp{\mat{D}}\mat{D}&-\trsp{\mat{D}}\mat{D}\\
-\trsp{\mat{D}}\mat{D}&\trsp{\mat{D}}\mat{D}
\end{bmatrix} 
\end{equation}
Note that the matrix $\mat{B}$ defined in \eqref{eq-z-gamma-B} is symmetric.
The gradient of the quadratic form $ Q(\vec{z}) $ is
\begin{equation}\label{GPSR-2}
\nabla_{\vec{z}} Q(\vec{z})=\vec{\gamma}+\mat{B}\vec{z}.
\end{equation}
Thus the iterative scheme based on the steepest-descent can be
designed by \eqref{GPSR-2}, viz.
\begin{equation}\label{GPSR-3}
\vec{z}_{k+1}=\vec{z}_k-\alpha_{k}\nabla Q(\vec{z}_k)
\end{equation}
where the optimal step size $ \alpha_{k} $  can be solved by a standard line-search process such as the Armijo rule \cite{bertsekas2016nonlinear}. Under appropriate convergence criteria, we can obtain the solution to \eqref{GPSR-1} according to \eqref{GPSR-3} by simple iterations.

The second GP method, known as the TNIPM method \cite{kim2007interior}, can be obtained by transforming the equation \eqref{QP_lambda} into the following constrained quadratic programming problem 
\begin{equation}\label{TNIPM-1}
\begin{aligned}
&\min_{\vec{r}\in \ES{R}{m}{1}} \frac{1}{2}\normp{\mat{D}\vec{r}-\vec{w}}{2}^2+\lambda\sum_{i=1}^{m}u_i,\\
&\st\ -u_i\le r_i\le u_i,\quad i=1,2,\cdots,m.
\end{aligned}
\end{equation}
By constructing the following \textit{logarithmic barrier} function for the bound constrain
\begin{equation}\label{TNIPM-2}
\Gamma(\vec{r},\vec{u})=-\sum_{i=1}^{m}\log(u_i+r_i)-\sum_{i=1}^{m}\log(u_i-r_i)
\end{equation}
over the domain 
\begin{equation}
\Omega = \set{(\vec{r},\vec{u})\in\mathbbm{R}^{m}\times\mathbbm{R}^{m}: \abs{r_i}\le u_i,i=1,2,\cdots,m} 
\end{equation}
and substituting \eqref{TNIPM-2} into \eqref{TNIPM-1}, we can find the unique optimal solution $ (\vec{r}^{*}(c),\vec{u}^{*}(c)) $ to the convex function \cite{kim2007interior}
\begin{equation}\label{TNIPM-3}
F_{c}(\vec{r},\vec{u})=c\cdot\left(\frac{1}{2}\normp{\mat{D}\vec{r}-\vec{w}}{2}^2+\lambda\sum_{i=1}^{m}u_i\right)+\Gamma(\vec{r},\vec{u})
\end{equation}
where $ c\in[0,+\infty) $. The optimal search direction $ \trsp{[\Delta\vec{r}, \Delta\vec{u}]}\in \Omega$ for \eqref{TNIPM-3} can be determined by 
\begin{equation}\label{TNIPM-4}
\nabla^{2}F_c(\vec{r},\vec{u})\cdot\begin{bmatrix}
\Delta \vec{r} \\ \Delta \vec{u}
\end{bmatrix}	 =-\nabla F_c(\vec{r},\vec{u})\in\ES{R}{2m}{1}
\end{equation}
with the Newton's method, in which the $\nabla^{2}F_c(\vec{r},\vec{u})$ in \eqref{TNIPM-4} is the Hessian matrix of $F_c(\vec{r},\vec{u})$. The optimal  step 
\begin{equation} \label{eq-step-tau}
\scrd{\tau}{opt} = \arg \min_{\tau\in \mathbb{R}} F_c(\vec{r} + \tau \Delta \vec{r}, \vec{u} + \tau \Delta\vec{u})
\end{equation}
specified by \eqref{eq-step-tau}
can be computed with a backtracking line search. In \cite{kim2007interior}, the search direction mentioned in \eqref{TNIPM-4} is accelerated by the \textit{preconditioned conjugate gradients} (PCG) algorithm \cite{nocedal1999numerical} to efficiently approximate the Hessian matrix.

\subsection{Homotopy Method}\label{sec:HPL1Opt}

The algorithm starts with an initial value $ \vec{r}^{(0)}=\vec{0} $ and operates iteratively, calculating the solutions $ \vec{r}^{(k)} $ at each step $ k=1,2,\cdots $. Throughout the computation, the active set
\begin{equation}\label{homotopy-1}
\mathscr{L}=\set{j\in\mathbb{N}:\abs{c^{(k)}_j}=\normp{\vec{c}^{(k)}}{\infty}=\lambda}
\end{equation}
remains constant, where $ \vec{c}^{(k)}=\trsp{\mat{D}}(\vec{w}-\mat{D}\vec{r}^{(k)}) $ is the vector of residual correlations \cite{ZhangXD2017Matrix}. Let
\begin{equation*}\label{homotopy-2}
(\cdot)_\lambda=(\cdot)[\mathscr{L}]
\end{equation*}
be the update direction on the sparse support, then $ \vec{d}^{(k)}_\lambda $ is the solution to the following linear system
\begin{equation}\label{homotopy-3}
\trsp{\mat{D}}_{\lambda}\mat{D}_{\lambda}\vec{d}^{(k)}_\lambda=\sign(\vec{c}^{(k)}_\lambda).
\end{equation}
Thus we can solve \eqref{homotopy-3} to obtain the $ \vec{d}^{(k)}_\lambda=\vec{d}^{(k)}[\mathscr{L}] $, which consists of the non-zero elements in $ \vec{d}^{(k)} $. Along the direction indicated by $ \vec{d}^{(k)} $, there are two scenarios when an update on $ \vec{r} $ may lead to a break-point \cite{asif2014homotopy}: 
inserting an element into the set $ \mathscr{L} $ or removing an element from the set $ \mathscr{L}$.
\blue{Let
	\begin{equation}\label{homotopy-4}
	\gamma_{+}^{(k)}=\min\limits_{i\in\mathscr{L}^{c}}\set{\min\left(\frac{\lambda-c^{(k)}_i}{1-\trsp{\mat{D}}_i\mat{D}_\lambda\vec{d}_\lambda^{(k)}},\frac{\lambda+c^{(k)}_i}{1+\trsp{\mat{D}}_i\mat{D}_\lambda\vec{d}_\lambda^{(k)}}\right)}_{+}
	\end{equation}}
and
\begin{equation}\label{homotopy-5}
\gamma_{-}^{(k)}=\min\limits_{i\in\mathscr{L}}\set{{-r^{(k)}_i}\big/{d^{(k)}_i}}_{+}
\end{equation}
\blue{where $ \min\set{\cdot}_{+}$} means that the minimum is taken over only positive arguments. The Homotopy algorithm proceeds to the next break-point and updates the sparse support set $ \mathscr{L} $ iteratively by the following  formula
\begin{equation}\label{homotopy-6}
\vec{r}^{(k+1)}=\vec{r}^{(k)}+\min\set{\gamma_{+}^{(k)},\gamma_{-}^{(k)}}\cdot\vec{d}^{(k)}, \quad k = 1, 2, \cdots
\end{equation}
The iteration terminates when $\normp{\vec{r}^{(k+1)}-\vec{r}^{(k)}}{2} \le \scrd{\epsilon}{re}\normp{\vec{r}^{(k)}}{2} $ for the given relative error $\scrd{\epsilon}{re}$ according to the Cauchy's convergence criteria.

\subsection{Iterative Shrinkage-Thresholding Method}\label{sec:ISTL1Opt}

The IST method \cite{daubechies2004iterative,wright2009sparse} are proposed to solve the \eqref{QP_lambda} as a special case of the following \textit{composite objective function}
\begin{equation}\label{IST-1}
\scrd{\vec{r}}{opt} = \arg \min_{\vec{r}\in \ES{R}{m}{1}} f(\vec{r})+\lambda g(\vec{r})
\end{equation}
where $ f(\vec{r})=\frac{1}{2}\normp{\mat{D}\vec{r}-\vec{w}}{2}^2 $ and $ g(\vec{r})=\normp{\vec{r}}{1} $. The updating rule of IST to minimize \eqref{IST-1} can be calculated by taking a second-order approximation of $ f $. Mathematically, we have \cite{wright2009sparse}
\begin{equation}\label{IST-2}
\begin{aligned}
\vec{r}^{(k+1)} = &\arg \min_{\vec{r}\in \ES{R}{m}{1}} \left\{f(\vec{r}^{(k)})+\trsp{(\vec{r}-\vec{r}^{(k)})}\cdot\nabla f(\vec{r}^{(k)})\right.\\
&\left.+\frac{1}{2}\normp{\vec{r}-\vec{r}^{(k)}}{2}^2\cdot\nabla^2f(\vec{r}^{(k)})+\lambda g(\vec{r})\right\}\\
\approx &\arg \min_{\vec{r}\in \ES{R}{m}{1}} \left\{\trsp{(\vec{r}-\vec{r}^{(k)})}\cdot\nabla f(\vec{r}^{(k)})\right.\\
&\left.+\frac{\alpha_{k}}{2}\normp{\vec{r}-\vec{r}^{(k)}}{2}^2+\lambda g(\vec{r})\right\}\\
= &\arg \min_{\vec{r}\in \ES{R}{m}{1}}\set{\frac{1}{2}\normp{\vec{r}-\vec{u}^{(k)}}{2}^2+\frac{\lambda}{\alpha_{k}} g(\vec{r})}
\end{aligned}
\end{equation}
where
\begin{equation}\label{IST-3}
\vec{u}^{(k)}=\vec{r}^{(k)}-\frac{1}{\alpha_{k}}\nabla f(\vec{r}^{(k)}), \quad \alpha_{k}\in\mathbbm{R}^{+}.
\end{equation}
Note that $ g(\vec{r}) $ is a separable function. As a result, a closed-form solution for $ \vec{r}^{(k+1)} $ can be obtained for each component
\begin{equation}\label{IST-4}
\begin{aligned}
r^{(k+1)}_i&=\arg \min_{r_i\in\mathbbm{R}} \set{\frac{(r_i-u^{(k)}_i)^2}{2}+\frac{\lambda}{\alpha_{k}}\abs{r_i}}\\
&=\textrm{soft}\left(u^{(k)}_i,\frac{\lambda}{\alpha_{k}}\right)
\end{aligned}
\end{equation}
where
\begin{equation}\label{IST-5}
\textrm{soft}(u,a)=\sign(u)\cdot\max\set{\abs{u}-a,0}
\end{equation}
is the \textit{soft-thresholding} or \textit{shrinkage} function \cite{donoho1995noising}. The coefficient $ \alpha_{k} $ used in \eqref{IST-2},\eqref{IST-3} and \eqref{IST-4} is adopted to approximate the Hessian matrix $ \nabla^2f $ and it can be calculated with the \textit{Barzilai-Borwein} spectral approach \cite{wright2009sparse}.

\subsection{Alternating Direction Method}\label{sec:ADML1Opt}

With the help of the auxiliary variable $ \vec{u}\in\ES{R}{(m-n)}{1} $, the problem \eqref{BPDN} can be converted to the following equivalent form \cite{yang2011alternating}
\begin{equation}\label{ADM-1}
\min\limits_{\vec{r}\in\ES{R}{m}{1},\vec{u}\in\ES{R}{(m-n)}{1}}\normp{\vec{r}}{1},\ \st\left\{\begin{aligned}
&\mat{D}\vec{r}-\vec{w}=\vec{u},\\
&\normp{\vec{u}}{2}\le\epsilon.
\end{aligned}\right.
\end{equation}
Moreover, we have an augmented Lagrangian subproblem of the form
\begin{equation}\label{ADM-2}
\begin{aligned}
&\min_{\vec{r}\in \ES{R}{m}{1},\vec{u}\in\ES{R}{(m-n)}{1}}&\norm{\vec{r}}_1-\trsp{\vec{y}}(\mat{D}\vec{r}+\vec{u}-\vec{w})+\\
&&\frac{\mu}{2}\normp{\mat{D}\vec{r}+\vec{u}-\vec{w}}{2}^2\\
&\st \ \normp{\vec{u}}{2}\le\epsilon&
\end{aligned}
\end{equation}
where $ \vec{y}\in\ES{R}{(m-n)}{1} $ is a multiplier and $ \mu>0 $ is a penalty parameter. Applying inexact alternating minimization to \eqref{ADM-2} yields the following iterative scheme \cite{yang2011alternating}
\begin{equation}\label{ADM-3}
\left\{\begin{aligned}
\vec{u}^{(k+1)}&=\mathscr{P}_{\vec{B}_\epsilon}\left(\frac{\vec{y}^{(k)}}{\mu}-\left(\mat{D}\vec{r}^{(k)}-\vec{w}\right)\right),\\
\vec{r}^{(k+1)}&=\textrm{soft}\left(\vec{r}^{(k)}-\tau\vec{g}^{(k)},\frac{\tau}{\mu}\right),\\
\vec{y}^{(k+1)}&=\vec{y}^{(k)}-\zeta\mu\left(\mat{D}\vec{r}^{(k+1)}+\vec{u}^{(k+1)}-\vec{w}\right).
\end{aligned}\right.
\end{equation}
where $ \tau > 0 $ is a proximal parameter, $ \zeta > 0 $ is a constant and
\begin{equation}\label{ADM-4}
\vec{g}^{(k)}=\trsp{\mat{D}}\left(\mat{D}\vec{r}^{(k)}+\vec{u}^{(k+1)}-\vec{w}-\frac{\vec{y}^{(k)}}{\mu}\right).
\end{equation}
Note that $ \mathscr{P}_{B(\epsilon)}: \ES{R}{(m-n)}{1}\to B(\epsilon) $  is the projection onto the closed ball 
$ B(\epsilon)= \set{\vec{d}\in\ES{R}{(m-n)}{1}:\normp{\vec{d}}{2}\le\epsilon}$ with radius $\epsilon$.
It should be pointed out that, when $ \epsilon=0 $, \eqref{ADM-4} can rewritten by
\begin{equation}\label{ADM-5}
\left\{\begin{aligned}
\vec{r}^{(k+1)}&=\textrm{soft}\left(\vec{r}^{(k)}-\tau\trsp{\mat{D}}\left(\mat{D}\vec{r}^{(k)}-\vec{w}-\frac{\vec{y}^{(k)}}{\mu}\right),\frac{\tau}{\mu}\right),\\
\vec{y}^{(k+1)}&=\vec{y}^{(k)}-\zeta\mu\left(\mat{D}\vec{r}^{(k+1)}-\vec{w}\right).
\end{aligned}\right.
\end{equation}
which induces a simple iterative algorithm for solving the problem \eqref{BP}.

\end{appendix}


\end{document}